\parskip=1ex plus 2pt minus 1pt
\baselineskip=17pt plus 4pt minus 2pt

\font\largest=cmr10 scaled\magstep2
\font\larger=cmr10 scaled\magstep1

\input amssym.def
\def\gg{\frak{g}}

\def\ad{\mathop{\rm ad}\nolimits}
\def\Aut{\mathop{\rm Aut}\nolimits}
\def\can{\mathop{\rm can}\nolimits}
\def\Coker{\mathop{\rm Coker}\nolimits}

\def\Cup{\mathop\cup\limits}
\def\Diff{\mathop{\rm Diff}\nolimits}
\def\Hitch{\mathop{\rm Hitch}\nolimits}
\def\Hom{\mathop{\rm Hom}\nolimits}

\def\Im{\mathop{\rm Im}\nolimits}
\def\Inv{\mathop{\rm Inv}\nolimits}
\def\Ker{\mathop{\rm Ker}\nolimits}
\def\Mor{\mathop{\rm Mor}\nolimits}
\def\Op{\mathop{\rm Op}\nolimits}
\def\Opscript{{\cal OP}}
\def\res{\mathop{\rm res}\limits}
\def\Spec{\mathop{\rm Spec}\nolimits}
\def\Supp{\mathop{\rm Supp}\nolimits}
\def\Sum{\sum\limits}

\def\limr#1{\smash{\mathop{\rm lim}\limits_{\displaystyle 
                   \rightarrow \atop \raise 1ex \hbox{$\scriptstyle #1$}}\,}}
\def\liml#1{\smash{\mathop{\rm lim}\limits_{\displaystyle 
                   \leftarrow \atop \raise 1ex \hbox{$\scriptstyle #1$}}\,}}

\def\maplongright#1{\smash{\mathop{\hbox to .45in {\rightarrowfill}}
                    \limits^{\lower 1ex \hbox{$\scriptstyle #1$}} }}
\def\mapright#1{\smash{\mathop{\longrightarrow}\limits^{#1}}}

\def\maplongleft#1{\smash{\mathop{\hbox to .45in {\leftarrowfill}}
                    \limits^{\lower 1ex \hbox{$\scriptstyle #1$}} }}
\def\mapdown#1{\Big\downarrow\rlap{$\vcenter{\hbox{$\scriptstyle#1$}}$}}
\def\maplongdown#1{\Bigg\downarrow\rlap{$\vcenter{\hbox{$\scriptstyle#1$}}$}}

\def\maptypedown#1#2{\llap{$\vcenter{\hbox{$\scriptstyle#1$}}$}
                  \Big\downarrow\rlap{$\vcenter{\hbox{$\scriptstyle#2$}}$}}
\def\maptwodown{\rlap{\Big\downarrow}{\lower .5ex \hbox{$\vee$}} }
\def\mapdoubledown{\rlap{\Big\downarrow}
                        {\raise .05ex \hbox{$\hskip.01em \big\downarrow$}} }
\mathchardef\tophook="317F
\def\maphookdown{\rlap{\smash{\raise 1.06ex \hbox{$\hskip .05em \tophook$}}}
                      {\smash{\lower .3ex \hbox{$\Big\downarrow$} }} }
\def\maplonghookdown{\rlap{\smash{\raise .5ex \hbox{$\hskip .05em \tophook$}}}
                      {\smash{\lower 2.3ex \hbox{$\Bigg\downarrow$} }} }

\def\mapvert{\Big\|}
\def\mapsubset{\hskip -0.7em \lower 0.5ex \hbox{
                 $\vcenter{\hbox{$\smash{\subset}$}} 
                  \rlap{$\vcenter{\hbox{$
                     {\smash{\hskip -0.6em \lower 0.285ex \hbox{---}}} 
                      \atop
                     {\smash{\hskip -0.6em \raise 0.285ex \hbox{---}}} 
                  $}}$} $ }}
\def\mapse#1{\longsearrow\rlap{$\vcenter{\hbox{$\hskip -1.0em\scriptstyle#1$}}$}}
\def\mapsw#1#2{\llap{$\vcenter{\hbox{$\hskip -.8em\scriptstyle#1$}}$}
               \longswarrow\rlap{$\vcenter{\hbox{$\hskip -0.5em\scriptstyle#2$}}$}}

\def\mapnw#1{\longnwarrow\rlap{$\vcenter{\hbox{$\hskip -1.0em\scriptstyle#1$}}$}}

\def\longsearrow{\mathop 
                  {\vcenter{ 
                    \hbox {$\backslash$ } \vskip -0.6ex 
                     \hbox {$\hskip .322em \backslash$ } \vskip -0.6ex 
                      \hbox {$\hskip .644em \backslash$ } \vskip -2.2ex 
                        \hbox {$\hskip .66em \triangleleft $ }  }} }

\def\longswarrow{\mathop 
                  {\vcenter{ 
                        \hbox {/} \vskip -0.6ex 
                       \hbox {\hskip -0.322em /} \vskip -0.6ex 
                      \hbox {\hskip -0.644em /} \vskip -2.1ex 
                     \hbox {$\hskip -0.66em \triangleright $ }  }} }

\def\longnwarrow{\mathop 
                  {\vcenter{ 
                        \hbox {$\triangleleft $} \vskip -2.1ex 
                      \hbox {\hskip -0.02em /}  \vskip -0.6ex 
                     \hbox {\hskip -0.342em /} \vskip -0.6ex 
                    \hbox {\hskip -0.664em /}   }} }

\def\righttypearrow#1#2{\smash{\mathop {\longrightarrow} 
     \limits^{#1} \limits_{ {\raise 1ex \hbox {$\scriptstyle#2$} } }  }}

\def\relbuilder#1#2{\smash {\mathop {\raise .5ex \hbox{$#1$}} 
          \limits_{ {\raise 2ex \hbox {$\scriptstyle#2$} } }  }}

\def\minimatrix#1#2#3#4{ \bigl( {{\scriptstyle#1} \atop {\scriptstyle#2}}
                                {{\scriptstyle#3} \atop {\scriptstyle#4}} \bigr) }

\def\sqr#1#2{{\vcenter{\vbox{\hrule height.#2pt
                \hbox{\vrule width.#2pt height#1pt \kern#1pt
                      \vrule width.#2pt}
                \hrule height.#2pt}}}}
\def\square{\mathchoice\sqr{4.9}2\sqr{4.9}2\sqr{2.1}2\sqr{1.5}2}

\def\boxtimes{\,{\mathbin{\rlap
                         {$\vcenter{\hbox{$\square$}}$}
                         {\hskip-.121545em\times} }\,}}

\def\uphat{{\smash{\lower .6ex \hbox{\^{}}}}}

\def\blackbox{\quad\hbox{\vrule height6pt width6pt}}

\def\chapter{\copyright}

\outer\def\beginsection#1 #2 \par{\vskip0pt plus .2\vsize\penalty-250
                          \vskip0pt plus-.2\vsize\bigskip\vskip\parskip
                          \message{#1 #2}\leftline
                                  {\indent\indent 
                                   \chapter \S \bf #1. \enspace #2.}
                          \nobreak\smallskip\noindent}

\outer\def\subsection#1{\medbreak\noindent {\bf \chapter.#1} \enspace}
\outer\def\nullsubsection{\medbreak\noindent}

\outer\def\beginproclaim#1{{\bf #1. \enspace} \begingroup \sl } 
\outer\def\endproclaim{\endgroup \par \medbreak} 

\outer\def\beginproeq{\endgroup} 
\outer\def\endproeq{\begingroup \sl } 

\outer\def\beginremarks#1{{\bf #1. \enspace}} 
\outer\def\endremarks{\par \medbreak} 

\outer\def\beginproof#1{\medbreak\noindent {\it Proof#1. \enspace} }
\outer\def\endproof{\blackbox \par \medbreak} 

\newcount\num
\def\eqnum#1#2{\global\advance\num by 1
             \eqno {(\chapter.\the\num {\rm #2})}
             \xdef#1{\the\num} }
\def\eqalignnum#1#2{\global\advance\num by 1
             {(\chapter.\the\num {\rm #2})}
             \xdef#1{\the\num} }

\def\label#1#2{{(\chapter.#1{\rm #2})}}

\newsymbol\twoheadrightarrow 1310


\widowpenalty=10000
\displaywidowpenalty=1000
\predisplaypenalty=5000
\postdisplaypenalty=100
\clubpenalty=10000

\magnification=1200


\centerline{\largest {Opers.}}
\smallskip
\centerline{\larger by A.\ Beilinson and V.\ Drinfeld}
\bigskip

In this article, which we wrote for ourselves in 1993,
there are no new results, we only rewrote some
definitions from [DS1] and [DS2] in a coordinate-free 
manner. We did this exercise because of the applications
of opers to the geometric Langalnds program.

We decided to publish this text as an e-print because it can
be used by students.

\beginsection 1 \hbox{$G\/$-opers} and \hbox{$\gg\/$-opers}

\subsection {1.1}
Let $G\/$ be a connected reductive group over ${\Bbb C}\/$ with a fixed Borel 
subgroup $B \subset G\/$ and a fixed Cartan subgroup $H \subset B\/$.  Let 
$\gg \supset {\frak b} \supset {\frak h}\/$ be the corresponding Lie algebras, 
$\Gamma \subset {\frak h}^*\/$ the set of simple roots with respect to 
${\frak b}\/$.  If $\alpha \in {\frak h}^*\/$ we set 
$\gg^\alpha = \{ x \in \gg | [a, x] = \alpha (a) x\/$ 
for all $a \in {\frak h}\/$'s.  There is a unique Lie algebra grading
$\gg = \relbuilder {\oplus}{k} \gg_k\/$ such that $\gg_0 = {\frak h}\/$, 
$\gg_1 = \relbuilder {\oplus}{\alpha \in \Gamma} \gg^\alpha\/$ and 
$\gg_{-1} = \relbuilder {\oplus}{\alpha \in \Gamma} \gg^{-\alpha}\/$.  The 
corresponding filtration $\gg^k = \relbuilder {\oplus}{r \geq k} 
\gg_r\/$ is \hbox{$B\/$-invariant}.

Fix a smooth algebraic curve $Y\/$ over ${\Bbb C}\/$.  Given a \hbox{$B\/$-bundle} 
${\cal F}\/$ on $Y\/$ we denote by ${\cal E_F}\/$ the algebroid of 
infinitesimal symmetries of ${\cal F}\/$ and by ${\cal E}^\gg_{\cal F}\/$ the 
algebroid of infinitesimal symmetries of the \hbox{$G\/$-bundle} corresponding to 
${\cal F}\/$.  We have a commutative diagram 
$$
\matrix{0 & \hskip -1.0em \strut \maplongright {} \hskip -1.0em & {\frak b}_{\cal F} & 
            \hskip -1.0em \maplongright {} \hskip -1.0em & {\cal E_F} & 
            \hskip -1.0em \maplongright {} \hskip -1.0em & 
            \Theta & \hskip -1.0em \maplongright {} \hskip -1.0em & 0 \cr
        & & \maphookdown & & \maphookdown & & \mapvert \cr
        0 & \hskip -1.0em \maplongright {} \hskip -1.0em & \gg_{\cal F} & 
            \hskip -1.0em \maplongright {} \hskip -1.0em & {\cal E}^\gg_{\cal F} & 
            \hskip -1.0em \maplongright {} \hskip -1.0em & \Theta & 
            \hskip -1.0em \maplongright {} \hskip -1.0em & 0 \cr
} \eqnum {\one} {}
$$
with exact rows where ${\frak b}_{\cal F}\/$ and $\gg_{\cal F}\/$ denote the 
\hbox{${\cal F}\/$-twists} of ${\frak b}\/$ and $\gg\/$ while $\Theta\/$ denotes the 
tangent bundle of $Y\/$.  We have the filtration $\gg^k_{\cal F} \subset 
\gg_{\cal F}\/$ and also the filtration ${\cal E}^k_{\cal F} \subset 
{\cal E}^\gg_{\cal F}, k \leq 0\/$, where ${\cal E}^k_{\cal F}\/$ 
denotes the preimage of $\gg^k_{\cal F} / {\frak b}_{\cal F} \subset 
\gg_{\cal F} / {\frak b}_{\cal F} = {\cal E}^\gg_{\cal F} / {\cal E_F}\/$.  
Notice that ${\cal E}^{-1}_{\cal F} / {\cal E_F}= \gg^{-1}_{\cal F} / 
\gg_{\cal F} = \relbuilder {\oplus}{\alpha \in \Gamma} \gg^{-\alpha}_{\cal F}\/$.  
Here $\gg^{-\alpha}_{\cal F}\/$ denotes the \hbox{${\cal F}\/$-twist} of the 
\hbox{$B\/$-module} $\gg^{-\alpha}\/$ (the action of $B\/$ on $\gg^{-\alpha}\/$ is 
defined to be the composition $B \rightarrow H \rightarrow \Aut \gg^{-\alpha}\/$).

\subsection {1.2}
\beginproclaim {Definition}
A \hbox{$G\/$-oper} on $Y\/$ is a \hbox{$B\/$-bundle} ${\cal F}\/$ on $Y\/$ with a
connection 
$\nabla : \Theta \rightarrow {\cal E}^\gg_{\cal F}\/$ such that 
$\nabla (\Theta) \subset {\cal E}^{-1}_{\cal F}\/$ and for every 
$\alpha \in \Gamma\/$ the composition 
$$
\Theta \righttypearrow{\nabla}{} 
{\cal E}^{-1}_{\cal F} \longrightarrow {\cal E}^{-1}_{\cal F} / {\cal E_F} 
\longrightarrow \gg^{-\alpha}_{\cal F}
$$ 
is an isomorphism.
\endproclaim

There is a natural notion of isomorphism of \hbox{$G\/$-opers}.  So 
\hbox{$G\/$-opers} on $Y\/$ form a groupoid. 

\nullsubsection
\beginremarks {Remark}
The term ``oper'' is motivated by the fact that for most of the classical $G\/$ 
one can interpret \hbox{$G\/$-opers} as differential operators between certain 
line bundles (see \chapter \S 2).  
\endremarks

\subsection {1.3}
Let $\gg\/$ be a semisimple Lie algebra over ${\Bbb C}\/$.

\nullsubsection
\beginproclaim {Definition}
A \hbox{$\gg\/$-oper} is a \hbox{$G\/$-oper} where $G\/$ is the group of inner 
automorphisms of $\gg\/$.
\endproclaim

\nullsubsection
\beginremarks {Remark}
The notion of \hbox{$\gg\/$-oper} (but not the word ``oper'') was introduced in [DS1] 
and [DS2] under the assumption that $\Theta\/$ is trivial and a trivialization 
of $\Theta\/$ is fixed.  The notion of \hbox{$sl(n)\/$-oper} for an arbitrary $Y\/$ 
was introduced in [Te1].  The notion of \hbox{$sl(2)\/$-oper} is equivalent to that of 
Sturm-Liouville operator or projective connection 
(see \chapter.2.6 and \chapter.2.7).  So \hbox{$sl(2)\/$-opers} have a long history 
(see the survey [Tyu]).
\endremarks

\nullsubsection
\beginproclaim {Proposition}
A \hbox{$\gg\/$-oper} does not have nontrivial automorphisms.
\endproclaim

\beginproof{}
The assertion is local, so we may assume that the \hbox{$B\/$-bundle} ${\cal F}\/$ 
corresponding to our \hbox{$\gg\/$-oper} is trivial and that there is a regular 
function $z\/$ on $Y\/$ such that $dz\/$ has no zeros.  The connection 
$\nabla\/$ on ${\cal F}\/$ is defined by an operator of the form 
$\partial_z + q\/$ where $q\/$ is a regular function on $Y\/$ with values 
in $\gg^{-1}\/$ such that for every $\alpha \in \Gamma\/$ the component 
$q_{-\alpha}\/$ of $q\/$ has no zeros.  An automorphism of 
$({\cal F}, \nabla)\/$ is a $B\/$-valued regular function $b\/$ on 
$Y\/$ such that 
$$
b^{-1} \left( \partial_z + q \right) b = \partial_z + q
\eqnum {\two} {}
$$
Let $q_k\/$ be the components of $q\/$ with respect to the grading 
$\gg = \relbuilder {\oplus}{k} \gg_k\/$.  Then \label{\two}{} implies that 
$b^{-1}_0 q_{-1} b_0 = q_{-1}\/$ where $b_0\/$ 
is the image of $b\/$ in $H\/$.  Since the center of $G\/$ is trivial 
$b_0 = 1\/$.  So we may consider the logarithm $\log b\/$.  
Let $u_k\/$ be the components of $\log b\/$ with 
respect to the grading $\gg = \relbuilder {\oplus}{k} \gg_k\/$.  
Obviously $u_0 = 0\/$.  Assume that $r > 0\/$ and $u_k = 0\/$ for 
$k < r\/$.  Then it follows from \label{\two}{} that $[q_{-1}, u_r] = 0\/$.  Notice 
that $q_{-1}\/$ is a member of an \hbox{$s l (2)\/$-triple} $\{ h, x, y \}\/$ with
$y = q_{-1}\/$ and $h \in {\frak h}\/$ being the sum of fundamental coweights.  Since 
$[h, u_r] = r u_r\/$ and $r > 0\/$ the equality $[y, u_r] = 0\/$ implies that 
$u_r = 0\/$.  So $u_k = 0\/$ for all $k\/$ and therefore $b = 1\/$.
\endproof

\nullsubsection
\beginproclaim {Corollary}
If $G\/$ is an arbitrary connected reductive group and $Y\/$ is connected then 
the group of automorphisms of any \hbox{$G\/$-oper} on $Y\/$ is equal to the center of 
$G\/$.
\blackbox
\endproclaim

\subsection{1.4}
Let $G\/$ and $\gg\/$ denote the same objects as in \chapter.1.1.  Denote by 
$G_{ss}\/$ the commutant of $G\/$ (i.e., the semisimple part of $G\/$) and by 
$\gg_{ss}\/$ the Lie algebra of $G_{ss}\/$.  Let us discuss the relation between 
\hbox{$G\/$-opers} and \hbox{$\gg_{ss}\/$-opers}.  Denote by $Z\/$ the center of $G\/$ and
by 
$\Omega\/$ the canonical bundle of $X\/$.  Of course to each \hbox{$G\/$-oper} there 
corresponds a \hbox{$(G / Z)\/$-oper}, i.e., a \hbox{$\gg_{ss}\/$-oper}.  On the other hand
if 
one fixes a square root of $\Omega\/$ (i.e., a line bundle ${\cal M}\/$ on $X\/$ with 
an isomorphism ${\cal M}^{\otimes 2} \righttypearrow {\sim}{} \Omega\/$) then there is 
a natural way to lift a \hbox{$\gg_{ss}\/$-oper} to a \hbox{$G_{ss}\/$-oper} and therefore to 
a \hbox{$G\/$-oper}.  Indeed, set $B_{\ad} = B / Z\/$, $H_{\ad} = H / Z\/$, 
$B_{ss} = B \cap G_{ss}\/$, $H_{ss} = H \cap G_{ss}\/$.  Let 
$({\cal F}_{\ad}, \nabla)\/$ be a \hbox{$\gg_{ss}\/$-oper}.  ${\cal F}_{\ad}\/$ is a 
\hbox{$B_{\ad}\/$-bundle} and lifting $({\cal F}_{\ad}, \nabla)\/$ to a \hbox{$G_{ss}\/$-oper} 
means lifting ${\cal F}_{\ad}\/$ to a \hbox{$B_{ss}\/$-bundle}.  
Denote by ${\cal F}^H_{\ad}\/$ the \hbox{$H_{\ad}\/$-bundle} corresponding to
${\cal F}_{\ad}\/$.  Lifting ${\cal F}_{\ad}\/$ to a \hbox{$B_{ss}\/$-bundle} means 
lifting ${\cal F}^H_{\ad}\/$ to a \hbox{$H_{ss}\/$-bundle}.
But for every $\alpha \in \Gamma\/$ we have a canonical isomorphism 
$\gg^{-\alpha}_{{\cal F}_{\ad}} \righttypearrow {\sim}{} \Theta\/$ 
(cf. \chapter.1.1) and therefore 
$\gg^{\alpha}_{{\cal F}_{\ad}} \righttypearrow {\sim}{} \Omega\/$.  Fix 
non-zero elements $y_\alpha \in \gg^{-\alpha}\/$.  Then the isomorphisms 
$\gg^{-\alpha}_{{\cal F}_{\ad}} \righttypearrow {\sim}{} \Theta\/$ induce an 
isomorphism between ${\cal F}^H_{\ad}\/$ and the direct image of the 
\hbox{${\Bbb G}_m\/$-torsor} on $Y\/$ corresponding to $\Omega\/$ with respect to the
homomorphism $\lambda : {\Bbb G}_m \rightarrow H_{\ad}\/$ such that 
$\lambda (t)\/$ acts on each $\gg^\alpha, \alpha \in \Gamma\/$, as 
multiplication by~$t\/$.  Notice that there is a unique homomorphism 
$\mu : {\Bbb G}_m \rightarrow H_{ss}\/$ such that the diagram 
$$
\matrix {{\Bbb G}_m & 
             \hskip -1.0em \vphantom{\Bigl(} \maplongright {t \mapsto t^2} \hskip -1.0em & 
             {\Bbb G}_m \cr
         \maptypedown {\mu}{} & & \maptypedown {}{\lambda} \cr
         H_{ss} & \hskip -1.0em \maplongright {} \hskip -1.0em & H_{\ad} \cr
}\eqnum {\three} {}
$$
is commutative (indeed, $\lambda\/$ corresponds to $\check{\rho}\/$, the sum 
of the fundamental coweights of $\gg_{ss}\/$, and since $2 \check{\rho}\/$ 
belongs to the coroot lattice $2 \check{\rho}\/$ corresponds to a homomorphism 
${\Bbb G}_m \rightarrow H_{ss})\/$.  So once $\Omega^{\otimes (1/2)}\/$ and 
$y_\alpha\/$ are fixed there is a natural way of lifting ${\cal F}^H_{\ad}\/$ 
to a \hbox{$H_{ss}\/$-bundle}: just take the direct image with respect to $\mu\/$ 
of the \hbox{${\Bbb G}_m\/$-torsor} corresponding to 
$\Omega^{\otimes (1/2)}\/$.  The dependence of the lifting on the choice of 
$y_\alpha\/$ is not very serious: the \hbox{$H_{ss}\/$-bundles} corresponding to 
$\{ y_\alpha \}\/$ and $\{ \widetilde{y}_\alpha \}\/$ 
$\bigl(\/$and to the same choice of $\Omega^{\otimes (1/2)} \bigr)\/$ 
are isomorphic but the isomorphism depends on the choice of $h \in H_{ss}\/$ 
such that $h y_\alpha h^{-1} = \widetilde{y}_\alpha\/$ for all $\alpha\/$.

So for a fixed $\Omega^{\otimes (1/2)}\/$ we have an equivalence between the 
category of \hbox{$G\/$-opers} and the category of pairs 
$({\cal E}_1, {\cal E}_2)\/$ where ${\cal E}_1\/$ is a 
\hbox{$\gg_{ss}\/$-oper} and ${\cal E}_2\/$ is a \hbox{$Z\/$-bundle} 
with a connection: given ${\cal E}_1\/$ one constructs a \hbox{$G_{ss}\/$-oper} 
as explained above and then one twists it by means of~${\cal E}_2\/$.  This 
equivalence depends on the choice of $\Omega^{\otimes (1/2)}\/$ in the 
following way.  Let ${\cal M}\/$ and $\widetilde{\cal M}\/$ be the square roots 
of~$\Omega\/$.  If the pair $({\cal E}_1, {\cal E}_2)\/$ corresponds to 
${\cal M}\/$ and a certain \hbox{$G\/$-oper} and 
$(\widetilde{\cal E}_1, \widetilde{\cal E}_2)\/$ corresponds to 
$\widetilde{\cal M}\/$ and the same \hbox{$G\/$-oper} then
$$
\widetilde{\cal E}_1 = {\cal E}_1, \; 
\widetilde{\cal E}_2 = {\cal E}_2 \otimes \nu_*
                   \bigl( \widetilde{\cal M} \otimes {\cal M}^{\otimes (-1)} \bigr)
\eqnum {\four} {}
$$
Here $\widetilde{\cal M} \otimes {\cal M}^{\otimes (-1)}\/$ is considered as a 
torsor over the group $\{ 1, -1 \}\/$ and 
$\nu_* \bigl( \widetilde{\cal M} \otimes {\cal M}^{\otimes (-1)} \bigr)\/$ 
is its image with respect to the homomorphism 
$\nu : \{ 1, -1 \} \rightarrow Z\/$ such that $\nu (-1) = \mu (-1)\/$ 
$\bigl(\/$the commutativity of \label{\three}{} implies that 
$\mu (-1) \in \Ker (H_{ss} \rightarrow H_{a d}) = H_{ss} \cap Z 
\bigr)\/$.

\beginsection 2 \hbox{$G\/$-opers} for classical $G\/$

We are going to show that for most of the classical groups $G\/$ a \hbox{$G\/$-oper}
can be considered as a differential operator between certain line bundles.

\subsection {2.1}
Let us start with the case $G = GL (n)\/$.  In this case a \hbox{$G\/$-oper} can be 
considered as a locally free \hbox{${\cal O}_Y\/$-module} sheaf of rank $n\/$ with a 
filtration ${\cal E = E}_n \supset {\cal E}_{n-1} \supset \ldots \supset 
{\cal E}_1 \supset {\cal E}_0 = 0\/$ and a connection $\nabla : {\cal E} 
\rightarrow {\cal E} \otimes \Omega\/$ such that

\itemitem{(i)}  the sheaves ${\cal E}_i / {\cal E}_{i-1}, 1 \leq i \leq n\/$,
                are invertible,

\itemitem{(ii)} {\rightskip = \parindent 
                for $i = 1, 2, \ldots, n-1 \; 
                \nabla ({\cal E}_i) \subset 
                {\cal E}_{i + 1} \otimes \Omega\/$ and $\nabla\/$ induces an 
                isomorphism    
                ${\cal E}_i / {\cal E}_{i-1} 
                \righttypearrow {\sim}{} ({\cal E}_{i+1} / {\cal E}_i) \otimes
                \Omega\/$.\par}

Denote by ${\cal D}\/$ the sheaf of differential operators on $Y\/$ and write 
${\cal O}\/$ instead of ${\cal O}_Y\/$.  Let ${\cal D}_i \subset {\cal D}\/$
be the sheaf of differential operators of order $\leq i\/$.  The connection 
$\nabla\/$ defines a \hbox{${\cal D}\/$-module} structure on ${\cal E}\/$.  So we have 
a natural morphism ${\cal D} \otimes_{\cal O} {\cal E}_1 \rightarrow {\cal E}\/$.
The condition (ii) implies that the morphism ${\cal D}_{n-1} \otimes_{\cal O} 
{\cal E}_1 \rightarrow {\cal E}\/$ is an isomorphism.  So considering the diagram
$$
\matrix {0 & \hskip -1.0em \mapright {} \hskip -1.0em & 
             {\cal D}_{n-1} \otimes_{\cal O} {\cal E}_1 & 
             \hskip -1.0em \mapright {} \hskip -1.0em & 
             {\cal D}_{n} \otimes_{\cal O} {\cal E}_1 &
             \hskip -1.0em \mapright {} \hskip -1.0em & 
             \Omega^{\otimes (-n)} \otimes_{\cal O} {\cal E}_1 &
             \hskip -1.0em \mapright {} \hskip -1.0em & 0    \cr
         & & \hfill \mapse {\sim} \hskip -2.0em & & 
             \hskip 0.5em \maplongdown {} \hfill    \cr
         & & & & {\cal E} \hfill       \cr
}
$$
we obtain a splitting $\Omega^{\otimes (-n)} \otimes_{\cal O} {\cal E}_1 
\rightarrow {\cal D}_{n} \otimes_{\cal O} {\cal E}_1\/$ and therefore a section
$L \in H^0 \bigl( Y, {\cal B} \otimes_{\cal O} {\cal D} \otimes_{\cal O} 
{\cal A}^{\otimes (-1)} \bigr)\/$ where ${\cal A} = {\cal E}^{\otimes (-1)}_1\/$,
${\cal B} = \Omega^{\otimes n} \otimes {\cal E}^{\otimes (-1)}_1\/$.  Consider
$L\/$ as a differential operator ${\cal A} \rightarrow {\cal B}\/$.  Its order
equals to $n\/$ and its principal symbol $\sigma \in H^0 \bigl( Y, 
{\cal A}^{\otimes (-1)} \otimes {\cal B} \otimes \Omega^{\otimes (-n)} \bigr)\/$
equals to~1.

Conversely, suppose we are given invertible sheaves ${\cal A, B}\/$ and a 
differential operator $L : {\cal A} \rightarrow {\cal B}\/$ of order $n\/$ 
whose principal symbol $\sigma \in H^0 \bigl( Y, {\cal A}^{\otimes (-1)} 
\otimes {\cal B} \otimes \Omega^{\otimes (-n)} \bigr)\/$ has no zeros (in
this situation $\sigma\/$ defines an isomorphism 
${\cal B} \righttypearrow{\sim}{} {\cal A} \otimes \Omega^{\otimes n}\/$ and
after identifying ${\cal B}\/$ with ${\cal A} \otimes \Omega^{\otimes n}\/$ 
we may say that $\sigma = 1\/$).  Then we construct a \hbox{$GL(n)\/$-oper} as follows.
We have a \hbox{${\cal D}\/$-module} morphism 
$r_L : {\cal D} \otimes_{\cal O} {\cal B}^{\otimes (-1)} \rightarrow 
{\cal D} \otimes_{\cal O} {\cal A}^{\otimes (-1)}\/$ defined by
$r_L (u) = u L\/$ where $L\/$ is
considered as a section of 
${\cal B} \otimes_{\cal O} {\cal D} \otimes_{\cal O} {\cal A}^{\otimes (-1)}\/$. 
Set ${\cal E} = \Coker  r_L\/$, ${\cal E}_i = 
\Im \bigl( {\cal D}_i \otimes_{\cal O} {\cal A}^{\otimes (-1)} \rightarrow 
{\cal E} \bigr)\/$ and consider the \hbox{${\cal D}\/$-module} ${\cal E}\/$ as an 
\hbox{${\cal O}\/$-module} with a connection $\nabla\/$.  Then 
$({\cal E}, \{ {\cal E}_i \}, \nabla)\/$ is an \hbox{$G L (n)\/$-oper}.

So we have constructed an {\sl equivalence between the category of 
\hbox{$GL (n)\/$-opers} and the category of triples $({\cal A, B}, L)\/$ where 
${\cal A}\/$ and ${\cal B}\/$ are invertible sheaves on $Y\/$ and 
$L : {\cal A} \rightarrow {\cal B}\/$ is a differential operator of order $n\/$
with non-vanishing principle symbol.\/}  Let us explain that by a morphism from 
$({\cal A}_1, {\cal B}_1, L_1)\/$ to $({\cal A}_2, {\cal B}_2, L_2)\/$ we
mean a pair of isomorphisms ${\cal A} \righttypearrow{\sim}{} {\cal A}_1\/$,
${\cal B} \righttypearrow{\sim}{} {\cal B}_1\/$ such that the diagram
$$
\matrix { {\cal A} & \hskip -1.0em \maplongright {L} \hskip -1.0em & 
                {\cal B}        \cr
          \mapdown {\wr}   & &          \mapdown {\wr} \cr
          {\cal A}_1 & \hskip -1.0em \maplongright {L_1} \hskip -1.0em & 
                {\cal B}_1      \cr
}
$$
is commutative.

We can make a slightly more precise statement.  Denote by $\Diff_n ({\cal A, B})\/$
the set of differential operators $L : {\cal A} \rightarrow {\cal B}\/$ of order
$n\/$ with non-vanishing principal symbol.  If $L \in \Diff_n ({\cal A, B})\/$ 
and $({\cal E}, \{ {\cal E}_i \}, \nabla)\/$ is the corresponding 
\hbox{$G L (n)\/$-oper} then we have canonical isomorphisms
$$
{\cal E}_i / {\cal E}_{i - 1} = {\cal E}_i \otimes \Omega^{\otimes (1-i)}
     = {\cal A}^{\otimes (-1)} \otimes \Omega^{\otimes (1-i)} = 
     {\cal B}^{\otimes (-1)} \otimes \Omega^{\otimes (n+1-i)}
\eqnum {\five} {}
$$
In particular ${\cal E}_1 = {\cal A}^{\otimes (-1)}\/$, 
${\cal E}/ {\cal E}_{n-1} = \Omega \otimes {\cal B}^{\otimes (-1)}\/$.
Denote by $\Op_n ({\cal A, B})\/$ the set of isomorphism classes of 
\hbox{$G L (n)\/$-opers} $({\cal E}, \{ {\cal E}_i \}, \nabla)\/$ with 
${\cal E}_1 = {\cal A}^{\otimes (-1)}\/$ and ${\cal E} / {\cal E}_{n-1} = 
\Omega \otimes {\cal B}^{\otimes (-1)}\/$ (i.e., isomorphisms 
${\cal E}_1 \righttypearrow{\sim}{} {\cal A}^{\otimes (-1)}\/$ and 
${\cal E} / {\cal E}_{n-1} = \Omega \otimes {\cal B}^{\otimes (-1)}\/$ are
fixed).  We have constructed a canonical bijection
$$
\Op_n ({\cal A, B}) \righttypearrow{\sim}{} \Diff_n ({\cal A, B})
\eqnum {\six} {}
$$

\subsection {2.2}
Recall that the transpose $L^t\/$ of a differential operator 
$L = \Sum_i f_i \partial_z^i\/$ is defined by 
$L^t = \Sum_i (- \partial_z)^i f_i\/$.  
In the global situation if $L\/$ is a differential operator 
${\cal A} \rightarrow {\cal B}\/$ where ${\cal A}\/$ and ${\cal B}\/$ are 
locally free \hbox{${\cal O}\/$-modules} then $L^t\/$ is a differential operator 
$\Omega \otimes {\cal B}^* \rightarrow \Omega \otimes {\cal A}^*\/$.  In 
particular transposition maps ${\cal D}\/$ to 
$\Omega \otimes_{\cal O} {\cal D} \otimes_{\cal O} \Omega^{\otimes (-1)}\/$,
i.e., the sheaf of differential operators $\Omega \rightarrow \Omega\/$.  The
induced morphism ${\cal A} \otimes_{\cal O} {\cal D} \otimes_{\cal O} {\cal B}
\rightarrow {\cal B} \otimes \Omega \otimes_{\cal O} {\cal D} \otimes_{\cal O} 
\Omega^{\otimes (-1)} \otimes {\cal A}\/$ will also be called transposition and 
denoted by $L \mapsto L^t\/$.

If ${\cal A}\/$ and ${\cal B}\/$ are invertible sheaves on $Y\/$ and 
$L \in \Diff_n ({\cal A, B})\/$ then 
$L^t \in \Diff_n \bigl( \Omega \otimes {\cal B}^{\otimes (-1)}, 
\Omega \otimes {\cal A}^{\otimes (-1)} \bigr)\/$.  On the other hand if 
$({\cal E}, \{ {\cal E}_i \}, \nabla)\/$ is a \hbox{$G L (n)\/$-oper} with 
${\cal E}_1 = {\cal A}^{\otimes (-1)}\/$ and ${\cal E} / {\cal E}_{n - 1} = 
\Omega \otimes {\cal B}^{\otimes (-1)}\/$ then we can consider the dual 
\hbox{$G L (n)\/$-oper} $({\cal E}^*, \{ ({\cal E}^*)_i \}, \nabla)\/$ where 
$({\cal E}^*, \nabla)\/$ is the dual of $({\cal E}, \nabla)\/$ and
$({\cal E}^*)_i = {\cal E}^\bot_{n-i} \subset {\cal E}^*\/$.  Since 
$({\cal E}^*)_1 = ({\cal E} / {\cal E}_{n - 1})^* = 
\Omega^{\otimes (-1)} \otimes {\cal B}\/$ and 
${\cal E}^* / ({\cal E}^*)_{n - 1} = ({\cal E}_1)^* = {\cal A}\/$ dualization 
defines a mapping 
$\Op_n ({\cal A, B}) \rightarrow \Op_n 
   \bigl( \Omega \otimes {\cal B}^{\otimes (-1)}, 
   \Omega \otimes {\cal A}^{\otimes (-1)} \bigr)\/$.

\nullsubsection
\beginproclaim {Proposition}
The diagram
$$
\matrix {\Diff_n ({\cal A, B}) & 
            \hskip -1.0em \vphantom{\Bigl(} \maplongright {L \mapsto -L^t} \hskip -1.0em & 
              \Diff_n 
                \bigl( \Omega \otimes {\cal B}^{\otimes (-1)}, 
                       \Omega \otimes {\cal A}^{\otimes (-1)} \bigr)   \cr
         \mapdown {\wr} & & \mapdown {\wr}    \cr
         \Op_n ({\cal A, B}) & \hskip -1.0em \maplongright {} \hskip -1.0em & 
              \Op_n 
                \bigl( \Omega \otimes {\cal B}^{\otimes (-1)}, 
                       \Omega \otimes {\cal A}^{\otimes (-1)} \bigr)   \cr
}
$$
is commutative.
\endproclaim

\beginproof{}
Let $({\cal E}, \{ {\cal E}_i \}, \nabla)\/$ and 
$(\widetilde{\cal E}, \{ \widetilde{\cal E}_i \}, \widetilde\nabla)\/$ be the
\hbox{$G L (n)\/$-opers} corresponding respectively to $L\/$ and $-L^t\/$.  We must
construct a horizontal non-degenerate pairing 
${\cal E} \otimes \widetilde{\cal E} \rightarrow {\cal O}\/$ such that 
$\widetilde{\cal E}_i = {\cal E}^\bot_{n-i}\/$ and the isomorphisms 
${\cal E} / {\cal E}_{n - 1} \righttypearrow{\sim}{} \widetilde{\cal E}^*_1\/$
and $\widetilde{\cal E} / \widetilde{\cal E}_{n - 1} \righttypearrow{\sim}{}
{\cal E}^*_1\/$ induced by the pairing are equal to the compositions 
${\cal E} / {\cal E}_{n - 1} \righttypearrow{\sim}{} 
\Omega \otimes {\cal B}^{\otimes (-1)} \righttypearrow{\sim}{} 
\widetilde{\cal E}^*_1\/$ and 
$\widetilde{\cal E} / \widetilde{\cal E}_{n - 1} \righttypearrow{\sim}{} 
{\cal A} \righttypearrow{\sim}{} {\cal E}^*_1\/$.  More precisely
the latter condition means the commutativity of the diagrams
$$
\matrix{ {\cal E} / {\cal E}_{n-1} = {\cal E}_1 \otimes \Omega^{\otimes (1-n)} 
             & \vphantom{\Bigl(} = & 
               {\cal A}^{\otimes (-1)} \otimes \Omega^{\otimes (1-n)}  \cr
         \hskip 7em \mapse {\sim}  & & \maplongdown {\varphi} \hskip 3em   \cr
         \hfill \widetilde{\cal E}^*_1 & = 
             & {\cal B}^{\otimes (-1)} \otimes \Omega \hfill \cr
}
\eqnum {\seven} {\thinspace a}
$$
$$
\matrix{ \widetilde{\cal E} / \widetilde{\cal E}_{n-1} = 
             \widetilde{\cal E}_1 \otimes \Omega^{\otimes (1-n)} 
             & \vphantom{\Bigl(} = & 
                {\cal B} \otimes \Omega^{\otimes (-n)}  \cr
         \hskip 7em \mapse {\sim} & & 
             \maplongdown {(-1)^{n-1} \psi} \hskip 4em \cr
         \hfill {\cal E}^*_1 & = & {\cal A} \hskip 4em  \cr
}
\global\advance\num by -1  
\eqnum {\seven} {\thinspace b}
$$
where the equalities denote the canonical isomorphisms \label{\four}{} while 
$\varphi\/$ and $\psi\/$ denote division by 
$\sigma \in H^0 \bigl( Y, {\cal A}^{\otimes (-1)} \otimes {\cal B} \otimes 
\Omega^{\otimes (-n)} \bigr)\/$, the principle symbol of $L\/$.  Let us explain 
that $(-1)^{n-1}\/$ appears in \label{\seven}{\thinspace b} because the principal symbol of 
$-L^t\/$ equals $(-1)^{n-1} \sigma\/$.

To construct the pairing we use the sheaf $\underline{\cal D}\/$ of 
pseudodifferential symbols $\bigl(\/$local sections of $\underline{\cal D}\/$ 
can be written as formal expressions 
$\Sum^k_{i = - \infty} f_i \partial_z^i\/$ or as 
$\Sum^k_{i = - \infty} \partial_z^i g_i\/$ where $f_i\/$
and $g_i\/$ are regular functions and $z\/$ is a local parameter on $Y \bigr)\/$.  
We also use the canonical morphism 
$\res : \underline{\cal D} \otimes_{\cal O} \Omega^{\otimes (-1)} \rightarrow
{\cal O}\/$ defined by $\bigl\{ \Sum_i f_i \partial_z^i \bigr\}
(dz)^{-1} \mapsto f_{-1}\/$.

Here is the construction.  We have ${\cal E} = \Coker r_L\/$ where 
$r_L : {\cal D} \otimes_{\cal O} {\cal B}^{\otimes (-1)} \rightarrow 
{\cal D} \otimes_{\cal O} {\cal A}^{\otimes (-1)}\/$ is defined by 
$r_L (u) = u L\/$, $L \in H^0 
\bigl( Y, {\cal B} \otimes_{\cal O} {\cal D} \otimes_{\cal O} 
{\cal A}^{\otimes (-1)} \bigr)\/$.  We also have 
$\widetilde{\cal E} = \Coker r_{L^t}\/$, 
$r_{L^t} : {\cal D} \otimes_{\cal O} {\cal A} \otimes \Omega^{\otimes (-1)} 
\rightarrow {\cal D} \otimes_{\cal O} {\cal B} \otimes \Omega^{\otimes (-1)}\/$, 
$L^t \in H^0 \bigl( Y, \Omega \otimes {\cal A}^{\otimes (-1)} \otimes_{\cal O} 
{\cal D} \otimes_{\cal O} {\cal B} \otimes \Omega^{\otimes (-1)} \bigr)\/$.
Define a pairing between ${\cal D} \otimes_{\cal O} {\cal A}^{\otimes (-1)}\/$
and ${\cal D} \otimes_{\cal O} {\cal B} \otimes \Omega^{\otimes (-1)}\/$ by
$$
(u, v) = \res \bigl( u L^{-1} v^t \bigr)
\eqnum {\eight} {}
$$
Here $u \in {\cal D} \otimes_{\cal O} {\cal A}^{\otimes (-1)}\/$ 
$\bigl(\/$i.e., $u\/$ is a local section of 
${\cal D} \otimes_{\cal O} {\cal A}^{\otimes (-1)} \bigr)\/$, 
$v \in {\cal D} \otimes_{\cal O} {\cal B} \otimes \Omega^{\otimes (-1)}\/$,
$v^t \in {\cal B} \otimes_{\cal O} {\cal D} \otimes_{\cal O} 
\Omega^{\otimes (-1)}\/$, $L^{-1} \in {\cal A} \otimes_{\cal O} 
\underline{\cal D} \otimes_{\cal O} {\cal B}^{\otimes (-1)}\/$ and therefore 
$u L^{-1} v^t \in \underline{\cal D} \otimes_{\cal O} \Omega^{\otimes (-1)}\/$.
The pairing \label{\eight}{} induces a pairing 
${\cal E} \otimes \widetilde{\cal E} \rightarrow {\cal O}\/$ with the required
properties.
\endproof

\subsection {2.3}
Let us discuss the case $G = S p (n)\/$.  A \hbox{$S p (n)\/$-oper} can be 
considered as a \hbox{$G L (n)\/$-oper} $({\cal E}, \{ {\cal E}_i \}, \nabla)\/$ 
(see \chapter.2.1) with a horizontal symplective form on ${\cal E}\/$ such that
${\cal E}^\bot_i = {\cal E}_{n-i}\/$.  We will construct an equivalence between 
the category of pairs $({\cal A}, L)\/$ where ${\cal A}\/$ is an invertible 
sheaf on $Y\/$ and $L\/$ is a differential operator 
${\cal A} \rightarrow \Omega \otimes {\cal A}^{\otimes (-1)}\/$ such that 
$L^t = L\/$ and the principle symbol $\sigma\/$ of $L\/$ does not vanish 
(in this case $\sigma\/$ induces an isomorphism 
$\Omega^{\otimes n} \otimes {\cal A} \righttypearrow{\sim}{} 
\Omega \otimes {\cal A}^{\otimes (-1)}\/$ and therefore 
$\bigl( {\cal A} \otimes \Omega^{\otimes (n/2)} \bigr)^{\otimes 2} 
\righttypearrow{\sim}{} \Omega\/$, so 
${\cal A} = \Omega^{\otimes (1/2)} \otimes \Omega^{\otimes (-n/2)}\/$ 
for some choice of $\Omega^{\otimes (1/2)}\/$.

If a pair $({\cal A}, L)\/$ is given then the \hbox{$G L (n)\/$-oper}
$({\cal E}, \{ {\cal E}_i \}, \nabla)\/$ constructed in \chapter.2.1 is a 
\hbox{$S p (n)\/$-oper}, the symplectic form on ${\cal E}\/$ being defined by
\label{\eight}{}.  It is skew-symmetric because $L^t = L\/$ and therefore 
$\res \bigl( u L^{-1} v^t \bigr) = 
\res \bigl( ( v L^{-1} u^t )^t \bigr) = 
-\res \bigl( v L^{-1} u^t \bigr)\/$.

Conversely, let $({\cal E}, \{ {\cal E}_i \}, \nabla)\/$ be a 
\hbox{$S p (n)\/$-oper} and $L : {\cal A} \rightarrow {\cal B}\/$ the 
differential operator constructed from $({\cal E}, \{ {\cal E}_i \}, \nabla)\/$ 
in \chapter.2.1.  According to \label{\five}{} 
${\cal A} = {\cal E}^{\otimes (-1)}_1\/$ and 
${\cal B} = \Omega \otimes ({\cal E} / {\cal E}_{n-1})^{\otimes (-1)}\/$.  The
pairing ${\cal E} / {\cal E}_{n-1} \otimes {\cal E}_1 \rightarrow {\cal O}\/$
induced by ${\cal E} \otimes {\cal E} \rightarrow {\cal O}\/$ induces an
isomorphism $({\cal E / E}_{n-1})^* \righttypearrow{\sim}{} {\cal E}_1\/$ and
therefore 
${\cal B} \righttypearrow{\sim}{} \Omega \otimes {\cal A}^{\otimes (-1)}\/$.
So we can consider $L\/$ as an operator 
${\cal A} \rightarrow \Omega \otimes {\cal A}^{\otimes (-1)}\/$.  To show that
$L^t = L\/$ notice that according to Proposition \chapter.2.2 the element of
$\Op_n \bigl( {\cal A}, \Omega \otimes {\cal A}^{\otimes (-1)} \bigr)\/$
corresponding to $-L^t \in 
\Diff_n \bigl( {\cal A}, \Omega \otimes {\cal A}^{\otimes (-1)} \bigr)\/$ is
the \hbox{$G L (n)\/$-oper} $({\cal E}^*, \{ ({\cal E}^*)_i \}, \nabla)\/$ with the
isomorphisms 
$\varphi : ({\cal E}^*)_1 \righttypearrow{\sim}{} {\cal A}^{\otimes (-1)}\/$
and $\psi : {\cal E}^* / ({\cal E}^*)_{n-1} \righttypearrow{\sim}{} {\cal A}\/$
defined as the compositions 
$({\cal E}^*)_1 \righttypearrow{\sim}{} ({\cal E/E}_{n-1})^* 
\righttypearrow{\sim}{} {\cal E}_1 \righttypearrow{\sim}{} 
{\cal A}^{\otimes (-1)}\/$ and ${\cal E}^* / ({\cal E}^*)_{n-1} 
\righttypearrow{\sim}{} ({\cal E}_1)^* \righttypearrow{\sim}{} {\cal A}\/$.  
Since $-L^t\/$ corresponds to 
$({\cal E}^*, \{ ({\cal E})^*_i \}, \nabla, \varphi, \psi)\/$ the operator 
$L^t\/$ corresponds to 
$({\cal E}^*, \{ ({\cal E})^*_i \}, \nabla, -\varphi, \psi)\/$.  The pairing
${\cal E} \otimes {\cal E} \rightarrow {\cal O}\/$ induces an isomorphism 
$({\cal E}, \{ {\cal E}_i \}, \nabla) \righttypearrow{\sim}{} 
({\cal E}^*, \{ ({\cal E}^*)_i \}, \nabla)\/$ such that diagrams
$$
\matrix{ {\cal E}_1 & \hskip -1.0em \maplongright {\sim} \hskip -1.0em & 
               ({\cal E}^*)_1  \cr
         \hfill \mapse {\sim} \hskip -1.4em & & \mapsw {\sim}{-\varphi} \hfill  \cr
         & {\cal A}^{\otimes (-1)}   \cr
}
\qquad
\matrix{ {\cal E}/{\cal E}_{n-1} & 
               \hskip -1.0em \maplongright {\sim} \hskip -1.0em & 
               {\cal E}^* / ({\cal E}^*)_{n-1}  \cr
         \hfill \mapse {\sim} \hskip -1.0em & & \mapsw {\sim}{\psi} \hfill  \cr
         & {\cal A}   \cr
}
$$
are commutative.  So $L^t = L\/$.

\subsection {2.4}
According to \chapter.1.4 \hbox{$s p (n)\/$-opers} may be considered as 
\hbox{$S p (n)\/$-opers} modulo twisting by \hbox{$Z\/$-torsors} where $Z\/$ is the center
of $S p (n)\/$.  According to \chapter.2.3 an \hbox{$S p (n)\/$-oper} is the same as 
a pair $({\cal M}, L)\/$ where ${\cal M}\/$ is a square root of $\Omega\/$ and 
$L\/$ is a differential operator 
${\cal M} \otimes \Omega^{\otimes (-n/2)} \rightarrow {\cal M} \otimes 
\Omega^{\otimes (n/2)}\/$ of order $n\/$ such that $L^t = L\/$ and
the principal symbol of $L\/$ equals~1.  Denote by $S ({\cal M})\/$ the set of 
all such $L\/$.  If ${\cal M} \otimes {\cal A}\/$ is another square root of 
$\Omega\/$ then ${\cal A}\/$ is a square root of ${\cal O}\/$ and the canonical 
connection on ${\cal A}\/$ induces a bijection 
$S ({\cal M}) \righttypearrow{\sim}{} S ({\cal M} \otimes {\cal A})\/$.  
Replacing ${\cal M}\/$ by ${\cal M} \otimes {\cal A}\/$ and 
$L \in S ({\cal M})\/$ by its image in $S ({\cal M} \otimes {\cal A})\/$ 
corresponds to twisting \hbox{$S p (n)\/$-opers} by \hbox{$Z\/$-torsors}.

So the set of \hbox{$s p (n)\/$-opers} can be identified with the set of differential
operators \hbox{$L : \Omega^{\otimes (1/2)} \otimes \Omega^{\otimes (-n/2)}\/$} 
$\rightarrow \Omega^{\otimes (1/2)} \otimes \Omega^{\otimes (n/2)}\/$ of
order $n\/$ such that $L^t = L\/$ and the principal symbol of $L\/$ equals 1
$\bigl(\/$the latter set does not depend on the choice of 
$\Omega^{\otimes (1/2)} \bigr)\/$.

\subsection {2.5}
Let $n \in {\Bbb N}\/$ be odd.  By an \hbox{$O (n)\/$-oper} we mean a \hbox{$G L
(n)\/$-oper}
$({\cal E}, \{ {\cal E}_i \}, \nabla)\/$ (see \chapter.2.1) with a horizontal 
non-degenerate symmetric bilinear form $B\/$ on ${\cal E}\/$ such that 
${\cal E}^\bot_i = {\cal E}_{n-i}\/$.  An \hbox{$S O (n)\/$-oper} in the sense of 
\chapter.1.2 can be considered as an \hbox{$O (n)\/$-oper} 
$({\cal E}, \{ {\cal E}_i \}, \nabla, B)\/$ with an isomorphism 
$\alpha : {\cal O} \righttypearrow {\sim}{} \det {\cal E}\/$ compatible with
the bilinear form on $\det {\cal E}\/$ induced by $B\/$.  If $L\/$ is a 
differential operator 
$\Omega^{\otimes (1-n)/2} \rightarrow \Omega^{\otimes (1+n)/2}\/$ of order 
$n\/$ such that $L^t = -L\/$ and the principal symbol of $L\/$ equals 1
then we construct an \hbox{$S O (n)\/$-oper} as follows.  Take the \hbox{$G L (n)\/$-oper}
$({\cal E}, \{ {\cal E}_i \}, \nabla)\/$ corresponding to $L\/$ according to
\chapter.2.1 and define the bilinear form $B\/$ on ${\cal E}\/$ by 
\label{\eight}{}.  Thus we obtain an \hbox{$O (n)\/$-oper}.  The canonical isomorphisms
${\cal E}_i / {\cal E}_{i -1} \righttypearrow{\sim}{} \Omega^{\otimes (n+1-2i)/2}\/$
$\bigr($cf. \label{\five}{}$\bigr)$ induce an isomorphism $\alpha : {\cal O} 
\righttypearrow{\sim}{} \det {\cal E}\/$.  To show that it is compatible with 
the bilinear form on $\det {\cal E}\/$ induced by $B\/$ we may assume that 
there is a regular function $z\/$ on $Y\/$ such that $dz\/$ has no zeros.  Then 
$\Omega\/$ is trivial and ${\cal E}\/$ can be identified with the sheaf of 
differential operators on $Y\/$ of order $\leq n-1\/$.  So there is a base 
$e_i = \partial_z^i 
\in H^0 (Y, {\cal E})\/$, $0 \leq i \leq n-1\/$.  The isomorphism 
$\alpha : {\cal O} \righttypearrow {\sim}{} \det {\cal E}\/$ is defined by
$1 \mapsto e_0 \wedge \ldots \wedge e_{n-1}\/$ and we have to show that the 
matrix $B_{ij} = B (e_i, e_j)\/$ has determinant~1.  Indeed, $B_{ij} = 0\/$ 
if $i + j < n - 1\/$ and $B_{ij} = (-1)^i\/$ if $i + j = n - 1\/$, so 
$\det (B_{ij}) = 1\/$.

Conversely, if $({\cal E}, \{ {\cal E}_i \}, \nabla, B, \alpha)\/$ is an
\hbox{$S O (n)\/$-oper} then we have canonical isomorphisms 
${\cal E}_{k+1} / {\cal E}_k \righttypearrow {\sim}{} \det {\cal E} 
\righttypearrow {\sim}{} {\cal O}\/$ where $k = (n - 1) / 2\/$.  By virtue of
\label{\five}{} we have canonical isomorphisms 
${\cal E}_1 \righttypearrow{\sim}{} \Omega^{\otimes (n-1)/2}\/$, 
${\cal E/ E}_{n-1} \righttypearrow{\sim}{} \Omega^{\otimes (1-n)/2}\/$.  So
$({\cal E}, \{ {\cal E}_i \}, \nabla)\/$ defines an operator $L \in 
\Diff_n \bigl( \Omega^{\otimes (1-n)/2}, \Omega^{\otimes (1+n)/2} \bigr)\/$
whose principal symbol equals~1.  The argument used in \chapter.2.3 shows
that $L^t = -L\/$.

So for an odd $n\/$ we have constructed a bijection between the set of 
\hbox{$S O (n)\/$-opers} and the set of differential operators 
$L : \Omega^{\otimes (1-n)/2} \rightarrow \Omega^{\otimes (1+n)/2}\/$ of
order $n\/$ such that $L^t = -L\/$ and the principal symbol of $L\/$ equals~1.

\subsection {2.6}
Since $sl(2) = sp(2) = o(3)\/$ one can interpret \hbox{$sl(2)\/$-opers} in two 
different ways.

1)\enspace According to \chapter.2.4 \hbox{$sp(2)\/$-opers} on $Y\/$ bijectively
correspond to Sturm-Liouville operators 
$L : \Omega^{\otimes (1/2)} \otimes \Omega^{\otimes (-1)} \rightarrow 
\Omega^{\otimes (1/2)} \otimes \Omega\/$ for a fixed choice of 
$\Omega^{\otimes (1/2)}\/$.  By a Sturm-Liouville operator we mean a 
differential operator $L\/$ of order 2 such that $L^t = L\/$ and the principal 
symbol of $L\/$ equals~1.  The set of Sturm-Liouville operators 
$\Omega^{\otimes (1/2)} \otimes \Omega^{\otimes (-1)} \rightarrow 
\Omega^{\otimes (1/2)} \otimes \Omega\/$ ``does not depend'' on the choice of
$\Omega^{\otimes (1/2)}\/$ (for different choices of $\Omega^{\otimes (1/2)}\/$
there is a canonical bijection between the corresponding sets of 
Sturm-Liouville operators).  Sturm-Liouville operators are also called 
{\it projective connections\/} on $Y\/$ (we remind the origin of this term
in \chapter.2.7).

2)\enspace According to \chapter.2.5 \hbox{$o (3)\/$-opers} (i.e., 
\hbox{$S O (3)\/$-opers})
bijectively correspond to differential operators 
$\widetilde{L} : \Omega^{\otimes (-1)} \rightarrow \Omega^{\otimes (2)}\/$ of
order 3 such that $\widetilde{L}^t = -L\/$ and the principal symbol of 
$\widetilde{L}\/$ equals~1.

We are going to describe the relation between the operators $L\/$ and 
$\widetilde{L}\/$ corresponding to the same \hbox{$sl (2)\/$-oper}.  We will 
also give a third interpretation of \hbox{$sl(2)\/$-opers}.

Let us introduce some notation.  If ${\cal L}\/$ and ${\cal M}\/$ are 
\hbox{${\cal O}\/$-modules} set
$$
{\cal L} \boxtimes {\cal M} = p^*_1 {\cal L} \otimes p^*_2 {\cal M}
\eqnum {\nine} {}
$$
where $p_1\/$ and $p_2\/$ are the projections $Y \times Y \rightarrow Y\/$.
If $a, b \in {\Bbb Z}\/$, $a \geq b\/$, then set
$$
{\cal P}_{a, b} ({\cal L, {\cal M}}) =
     ({\cal L} \boxtimes {\cal M}) (a \Delta) / ({\cal L} 
             \boxtimes {\cal M}) (b \Delta)
\eqnum {\ten} {}
$$
where $\Delta \subset Y \times Y\/$ is the diagonal.

If ${\cal E}\/$ and ${\cal F}\/$ are locally free \hbox{${\cal O}\/$-modules} denote
by $\Diff_{\leq n} ({\cal E, F})\/$ the space of differential operators 
${\cal E} \rightarrow {\cal F}\/$ of order $\leq n\/$.  There is a well known 
canonical isomorphism
$$
\Diff_{\leq n} ({\cal E, F}) \righttypearrow{\sim}{} H^0 
    \bigl( Y \times Y, {\cal P}_{n + 1, 0} 
    ( \Omega \otimes {\cal E}^*, {\cal F} ) \bigr)
\eqnum {\eleven} {}
$$
characterized by the following property:  if $K(x, y)\/$ is a section of 
${\cal P}_{n+1, 0} \bigl( \Omega \otimes {\cal E}^*, {\cal F} \bigr)\/$ then the 
corresponding $L \in \Diff_{\leq n} ({\cal E, F})\/$ is defined by 
$$
(L s) (y) = n! \res_{x = y} \bigl( K (x,y) s (x) \bigr)
\eqnum {\twelve} {}
$$
(i.e., $L\/$ is the ``integral operator'' corresponding to the ``kernel''
$n! K\/$).  Replacing $K (x,y)\/$ by $K (y,x)\/$ is equivalent to replacing 
$L : {\cal E} \rightarrow {\cal F}\/$ by 
$-L^t : \Omega \otimes {\cal F}^* \rightarrow \Omega \otimes {\cal E}^*\/$.
Notice that due to the factor $n!\/$ in \label{\ten}{} the diagram
$$
\matrix {\hfill \Diff_{\leq n} ({\cal E, F}) & 
             \hskip -1.0em \maplongright{\sim} \hskip -1.0em & 
             H^0 \bigl( Y \times Y, {\cal P}_{n + 1, 0} 
                 ( \Omega \otimes {\cal E}^*, {\cal F} ) \bigr)   \cr
         \hskip 2.0em \maptypedown{\sigma}{} & & \mapdown{}   \cr
         H^0 \bigl( Y, {\cal E}^* \otimes \Omega^{\otimes (-n)} \otimes 
                      {\cal F} \bigr) & \hskip -1.0em 
             \maplongleft{\scriptstyle\varphi \atop \scriptstyle\sim} 
                      \hskip -1.0em &
             H^0 \bigl( Y \times Y, {\cal P}_{n + 1, n} 
                 ( \Omega \otimes {\cal E}^*, {\cal F} ) \bigr)   \cr
}
\eqnum {\thirteen} {}
$$
is commutative where $\sigma\/$ associates to 
$L \in \Diff_{\leq n} ({\cal E,F})\/$ its principal symbol and $\varphi\/$ is
induced by the isomorphism ${\cal O}_{Y \times Y} \bigl( (n+1) \Delta \bigr) / 
{\cal O}_{Y \times Y} (n \Delta) \righttypearrow{\sim}{} 
\Omega^{\otimes (-n-1)}\/$ defined by\break
$(x - y)^{-n-1} f(x,y) \mapsto f(x,x) (dx)^{-n-1}\/$ (here we consider 
${\cal O}_{Y \times Y} \bigl( (n+1) \Delta \bigr) / 
{\cal O}_{Y \times Y} (n \Delta) \/$ as a sheaf on $\Delta = Y\/$).

Now we can describe the relation between the Sturm-Liouville operator $L\/$ 
and the operator $\widetilde{L}\/$ of order 3 corresponding to the same 
\hbox{$sl (2)\/$-oper}. As explained above $L\/$ corresponds to a skew-symmetric 
section of 
${\cal P}_{3,0} \bigl( \Omega^{\otimes (3/2)}, \Omega^{\otimes (3/2)} \bigr)\/$, 
which can be uniquely lifted to a skew-symmetric 
$K \in H^0 \bigl( Y \times Y, {\cal P}_{3, -1} 
    ( \Omega^{\otimes (3/2)}, \Omega^{\otimes (3/2)} ) \bigr)\/$.
On the other hand $\widetilde{L}\/$ corresponds to a symmetric 
$\widetilde{K} \in H^0 \bigl( Y \times Y, {\cal P}_{4, 0} 
    ( \Omega^{\otimes 2}, \Omega^{\otimes 2} ) \bigr)\/$.  Both
$K\/$ and $\widetilde{K}\/$ have standard principal parts:
$K (x,y) = (x - y)^{-3} (dx)^{3/2} (dy)^{3/2} + \ldots\/$,
$\widetilde{K} (x,y) = (x - y)^{-4} (dx)^{2} (dy)^{2} + \ldots\/$.

\nullsubsection
\beginproclaim {Proposition}
$\widetilde{K} = K^{4/3}\/$.
\endproclaim

\beginproof{}
Suppose that in terms of a local parameter $z\/$ the operator $L\/$ is written
as $\partial_z^2 + u\/$ 
$\bigl( \/$i.e., $L \bigl( f (dz)^{-1/2} \bigr) = 
\bigl(\partial^2_z f + u f \bigr) (dz)^{3/2} \bigl)\/$.  The corresponding 
\hbox{$s l (2)\/$-oper} is defined by the connection
$$
\partial_z + \left( \matrix{ 0  & -u   \cr
                                 1  &  0   \cr} \right)
\eqnum {\fourteen}{}
$$
By $o (3)\/$ we understand the Lie algebra of the automorphism group of the quadratic form 
corresponding to the matrix
$$
\left( \matrix{ 0  &  0  &  1   \cr
                0  & -1  &  0   \cr
                1  &  0  &  0   \cr} \right)
$$
The image of \label{\fourteen}{} under a suitable isomorphism 
$s l (2) \righttypearrow{\sim}{} o (3)\/$ equals to
$$
\partial_z + \left( \matrix{ 0  & -2u  &  0   \cr
                                 1  &  0   & -2u  \cr
                                 0  &  1   &  0   \cr} \right)
\eqnum {\fifteen}{}
$$
The \hbox{$o (3)\/$-oper} defined by \label{\fifteen}{} corresponds to 
$\widetilde{L} = \partial_z^3 + 2 u \cdot \partial_z + 2 \partial_z 
\cdot u\/$.  Recall that $2 K\/$ 
and $6 \widetilde{K}\/$ are the ``kernels'' of $L\/$ and $\widetilde{L}\/$.  So
$$
\eqalignno {
           K & = \left\{ {{1}\over{(z_1 - z_2)^3}} + 
                     {{u (z_1) + u (z_2)}\over{4 (z_1 - z_2)}} + \ldots
                 \right\}
                 (dz_1)^{3/2} (dz_2)^{3/2} &
                 \eqalignnum {\sixteen}{}       \cr   \noalign {\vskip 1.0ex}
          \widetilde{K} & = \left\{ {{1}\over{(z_1 - z_2)^4}} + 
                               {{u (z_1) + u (z_2)}\over{3 (z_1 - z_2)^2}} 
                               + \ldots
                            \right\}
                            (dz_1)^{2} (dz_2)^{2}   &
                            \eqalignnum {\seventeen}{}   \cr
        }
$$
Therefore $\widetilde{K} = K^{4/3}\/$.
\endproof

In \chapter.2.7 we will sketch a non-computational proof of the formula 
$\widetilde{K} = K^{4/3}\/$.

Now we can give the third interpretation of \hbox{$s l (2)\/$-opers}: they are 
in bijective correspondence
with symmetric ${\cal K} \in H^0 
\bigl( Y \times Y, {\cal P}_{2, -2} ( \Omega, \Omega ) \bigr)\/$
with the standard principle part 
${\cal K} (x, y) = (x - y)^{-2} dx \; dy + \ldots\/$.
The bijection is defined by 
$$
{\cal K} = \widetilde{K}^{1/2} = K^{2/3}
\eqnum {\eighteen}{}
$$
where $2 K\/$ is the ``kernel'' of the Sturm-Liouville operator 
$L : \Omega^{\otimes (-1/2)} \rightarrow \Omega^{\otimes (3/2)}\/$ 
corresponding to the \hbox{$s l (2)\/$-oper} and $6 \widetilde{K}\/$ is the ``kernel''
of the corresponding differential operator 
$\widetilde{L} : \Omega^{\otimes (-1)} \rightarrow \Omega^{\otimes 2}\/$ of
order 3.

\subsection {2.7}
In this section we remind various facts about \hbox{$s l (2)\/$-opers} while in 
\chapter.2.8 and \chapter.2.9 we give an interpretation of \hbox{$S L (n)\/$-opers}
and \hbox{$S O (2k)\/$-opers}.  The reader can skip \chapter.2.7--\chapter.2.9 and 
pass directly to \chapter\S 3. 

In the previous sections we considered {\it algebraic\/} opers.  If $Y\/$ is
an analytic curve over ${\Bbb C}\/$ and the words ``bundle'' and ``connection''
in Definition \chapter.1.2 are understood in the analytic sense then one obtains
the notions of {\it analytic\/} \hbox{$G${\it -oper\/}} and {\it analytic\/} 
\hbox{$\gg${\it -oper\/}}.  An analytic \hbox{$s l (2)\/$-oper} on $Y\/$ is the same as a 
{\it locally-projective structure\/} on $Y\/$.  Indeed, an \hbox{$s l (2)\/$-oper} 
(i.e., a \hbox{$P G L (2)\/$-oper}) can be considered as a quadruple 
$(E, \pi, s, \nabla)\/$ where $\pi : E \rightarrow Y\/$ is a locally trivial
bundle with ${\Bbb P}^1\/$ as a fiber, $\nabla\/$ is a connection on it and 
$s : Y \rightarrow E\/$ is a section such that the covariant differential
of $s\/$ does not vanish.  In the analytic situation this is the same as 
locally-projective structure.  This explains the term ``projective connection''
used as a synonym for ``Sturm-Liouville operator'' or \hbox{``$s l (2)\/$-oper.''}
In terms of the Sturm-Liouville operator $L\/$ ``good'' local coordinates on 
$Y\/$ are defined as $\psi_1 / \psi_2\/$ where $\psi_1\/$ and $\psi_2\/$ are
local solutions to the equation $L \psi = 0\/$; two good local coordinates 
$z\/$ and $\zeta\/$ are related by a projective transformation 
$z = (a \zeta + b) / (c \zeta + d)\/$.

In the analytic situation the relation between the differential operators 
$L : \Omega^{\otimes (-1/2)} \rightarrow \Omega^{\otimes (3/2)}\/$ and 
$\widetilde{L} : \Omega^{\otimes (-1)} \rightarrow \Omega^{\otimes 2}\/$
can be expressed in the following way:  the sheaf $\Ker \widetilde{L}\/$ 
is generated by the products $\psi_1 \psi_2\/$ where $\psi_1\/$ and $\psi_2\/$
are local sections of $\Ker L\/$.

Notice that $\Omega^{\otimes (-1)} = \Theta\/$ is the tangent sheaf and
therefore has a Lie algebra structure.  In fact 
$\Ker \widetilde{L} \subset \Theta\/$ is a Lie subalgebra which is nothing but 
the sheaf of vector fields preserving $L\/$ (i.e., preserving the 
locally-projective structure corresponding to $L\/$).  This is an immediate 
consequence of the following expression for $\widetilde{L}\/$ in terms of
$L\/$ which holds in the algebraic situation as well as in the analytic
one:
$$
\widetilde{L} (v) = 2 [v, L]
\eqnum {\nineteen}{}
$$
Here $v \in \Theta = \Omega^{\otimes (-1)}\/$ (i.e., $v\/$ is a local 
section of $\Theta\/$) and $[v, L]\/$ is the Lie derivative of $L\/$ 
with respect to $v\/$ $\bigl( \Theta\/$ acts on $\Omega^{\otimes (-1)}\/$, 
$\Omega^{\otimes (3/2)}\/$, and therefore on the sheaf of differential 
operators $\Omega^{\otimes (-1/2)} \rightarrow \Omega^{\otimes (3/2)}\/$; 
in fact the differential operator $[v, L]\/$ has order 0, i.e., 
$[v, L] \in \Omega^{\otimes 2} \bigl)\/$.

\nullsubsection
\beginremarks {Remark}
Proposition \chapter.2.6 can be proved without computation.  Indeed, the formula
$L \mapsto \widetilde{K} - K^{4/3}\/$ defines a mapping $\varphi\/$ : 
$\{${\it projective connections\/}$\}$ $\longrightarrow$ 
$\{${\it quadratic differentials\/}$\}$.
To show that for any $L\/$ the quadratic differential $\omega = \varphi (L)\/$ 
vanishes at any $y \in Y\/$ introduce a ``good'' local parameter $\zeta\/$ at 
$y\/$, i.e., $\zeta = \psi_1 / \psi_2\/$ where $\psi_1\/$ and $\psi_2\/$ are
local solutions to $L \psi = 0\/$ such that $\psi_1 (y) = 0\/$.  Then 
$L \bigl( f \cdot (d \zeta)^{-1/2}\bigr) = f'' (d \zeta)^{3/2}\/$ and $L\/$ is
invariant under $\zeta \mapsto a \zeta\/$.  So $\omega\/$ is also invariant and 
therefore $\omega (y) = 0\/$.  This proof is in fact purely
algebraic (one can consider formal solutions to $L \psi = 0\/$ instead of
analytic ones).
\endremarks

\subsection {2.8}
An \hbox{$S L (n)\/$-oper} can be considered as a \hbox{$G L (n)\/$-oper} 
$({\cal E}, \{ {\cal E}_i \}, \nabla)\/$ with a fixed horizontal isomorphism
$$
{\cal O} \righttypearrow {\sim}{} \det {\cal E}
\eqnum {\twenty}{}
$$
As explained in \chapter.2.1 $({\cal E}, \{ {\cal E}_i \}, \nabla)\/$ defines 
a differential operator 
$L : {\cal A} \rightarrow {\cal A} \otimes \Omega^{\otimes n}\/$ or order
$n\/$ with principal symbol 1 where ${\cal A} = {\cal E}^{\otimes (-1)}_1\/$.
The isomorphisms \label{\five}{} and \label{\twenty}{} induce an isomorphism 
$$
{\cal A}^{\otimes n} \otimes \Omega^{\otimes n (n-1)/2} 
        \righttypearrow {\sim}{} {\cal O}
\eqnum {\twentyone}{}
$$
$L^t\/$ is an operator 
${\cal A} \otimes {\cal M} \rightarrow {\cal A} \otimes \Omega^{\otimes n} 
\otimes {\cal M}\/$ where 
${\cal M} = \Omega^{\otimes (1-n)} \otimes {\cal A}^{\otimes -2}\/$.  The 
isomorphism \label{\twentyone}{} induces an isomorphism 
${\cal M}^{\otimes n} \righttypearrow {\sim}{} {\cal O}\/$ and therefore a connection
on ${\cal M}\/$.  So $L^t\/$ induces a differential operator 
${\cal A} \rightarrow {\cal A} \otimes \Omega^{\otimes n}\/$ which we also 
denote by $L^t\/$.

\nullsubsection
\beginproclaim {Proposition}
The order of $L + (-1)^{n+1} L^t\/$ is less than $n - 1\/$ (i.e., the 
subprincipal symbol of $L\/$ equals 0).  Thus one obtains an equivalence
between the category of \hbox{$S L (n)\/$-opers} and the category of pairs 
$({\cal A}, L)\/$ where ${\cal A}\/$ is an invertible sheaf on $Y\/$ with
an isomorphism \label{\twentyone}{} and 
$L : {\cal A} \rightarrow {\cal A} \otimes \Omega^{\otimes n}\/$ is a 
differential operator of order $n\/$ with principal symbol~1 and subprincipal
symbol~0.
\endproclaim

The proof is left to the reader.  It is based on the following local
observation.  Assume that there is a regular function $z\/$ on $Y\/$ such that 
$dz\/$ has no zeros.  Let $({\cal E}, \{ {\cal E}_i \}, \nabla)\/$ be the 
\hbox{$G L (n)\/$-oper} corresponding to the operator 
$L : {\cal O} \rightarrow {\cal O}\/$, 
$L = \partial_z^n + \Sum^{n-1}_{i=0} f_i \partial_z^i\/$.  Identify 
${\cal E}\/$ with the sheaf of differential operators on $Y\/$ of order
$\leq n - 1\/$.  Then the matrix $q\/$ of the connection $\nabla\/$ in the
basis $e_i = \partial_z^i \in H^0 (Y, {\cal E})\/$, $0 \leq i \leq n-1\/$, 
equals 
$$
- \left( \matrix{ 0 \ldots & 0      &  f_0     \cr
                  \ldots   & \ldots & \ldots   \cr
                  0 \ldots & 0      & f_{n-1}  \cr
                } \right)
$$
and so $T r q = - f_{n-1}\/$.

Clearly there is a natural bijection between the set of \hbox{$s l (n)\/$-opers} and 
the set of differential operators 
$L : \Omega^{\otimes (1-n)/2} \rightarrow \Omega^{\otimes (1+n)/2}\/$ of
order $n\/$ with principal symbol 1 and subprincipal symbol 0 (the latter 
set does not depend on the choice of $\Omega^{\otimes (1/2)}\/$).

\subsection {2.9}
Let us describe \hbox{$S O (n)\/$-opers} and \hbox{$s o (n)\/$-opers} for $n = 2 k\/$,
$k \in {\Bbb Z}\/$, $k \geq 2\/$.  By an $S O (n)\/$ vector bundle we mean a 
locally free \hbox{${\cal O}\/$-module} ${\cal F}\/$ of rank $n\/$ with a 
non-degenerate symmetric bilinear form $B\/$ on it and an isomorphism 
$\alpha : {\cal O} \righttypearrow{\sim}{} \det {\cal F}\/$ compatible with 
the bilinear form $\overline{B}\/$ on $\det {\cal F}\/$ induced by $B\/$ 
$\bigl( \/$if $Y\/$ were a curve over a field which is not algebraically closed 
then it would probably be more reasonable to require $\alpha\/$ to be 
compatible with $(-1)^k \overline{B}\/$ because the determinant of the 
hyperbolic symmetric bilinear form on ${\Bbb Z}^{2k}\/$ equals 
$(-1)^k \bigr)\/$.  An \hbox{$S O (n)\/$-oper} is an $S O (n)\/$ vector bundle 
${\cal F}\/$ with a filtration 
${\cal F} = {\cal F}_{n-1} \supset {\cal F}_{n-2} \supset \ldots \supset
{\cal F}_1 \supset {\cal F}_0 = 0\/$ and a connection 
$\nabla : {\cal F} \rightarrow {\cal F} \otimes \Omega\/$ compatible with 
$B\/$ such that

\itemitem{1)} {\rightskip = \parindent 
              the sheaves ${\cal F}_i / {\cal F}_{i-1}\/$ are invertible for
              $1 \leq i \leq n - 1\/$, $i \neq k\/$ while 
              ${\cal F}_k / {\cal F}_{k - 1}\/$ is locally free of rank 2; \par}

\itemitem{2)} ${\cal F}^\bot_i = {\cal F}_{n-1-i}\/$;

\itemitem{3)} $\nabla ({\cal F}_i) \subset {\cal F}_{i+1} \otimes \Omega\/$
              for $1 \leq i \leq n - 2\/$;

\itemitem{4)} {\rightskip = \parindent 
              $\nabla\/$ induces an isomorphism 
              ${\cal F}_i / {\cal F}_{i - 1} \righttypearrow{\sim}{} 
              \bigl( {\cal F}_{i+1} / {\cal F}_i \bigr) \otimes \Omega\/$ 
              \enspace for $1 \leq i \leq n - 2\/$, $i \neq k-1\/$, $k\/$; 
              \par}

\itemitem{5)} {\rightskip = \parindent 
              the composition of the morphisms 
              $$
              {\cal F}_{k-1} / {\cal F}_{k - 2} \rightarrow
              \bigl( {\cal F}_k / {\cal F}_{k - 1} \bigr) \otimes \Omega
              \rightarrow \hbox{$\bigl( {\cal F}_{k+1} / {\cal F}_k \bigr) 
              \otimes \Omega^{\otimes 2} $}
              $$ 
              induced by $\nabla\/$ is an
              isomorphism. \par}

Let $L : \Omega^{\otimes (1-k)} \rightarrow \Omega^{\otimes k}\/$ be a 
differential operator of  order $2 k - 1\/$ such that $L^t = -L\/$ and the 
principal symbol of $L\/$ equals~1.  Starting from $L\/$ and a section 
$f \in H^0 \bigl( Y, \Omega^{\otimes k} \bigr)\/$ one can construct an 
\hbox{$S O (n)\/$-oper} in the following way.  Let 
$({\cal E}, \{ {\cal E}_i \}, \nabla_{\cal E})\/$ be the 
\hbox{$S O (n-1)\/$-oper} corresponding to $L\/$ (see \chapter.2.5).  According 
to \label{\five}{} we have canonical isomorphisms 
${\cal E}_i / {\cal E}_{i-1} \righttypearrow{\sim}{} \Omega^{\otimes (k-i)}\/$ 
and in particular ${\cal E}_1 \righttypearrow{\sim}{} \Omega^{\otimes (k-1)}\/$.
Denote by ${\cal F}\/$ the orthogonal direct sum ${\cal E} \oplus {\cal O}\/$, 
the bilinear form on ${\cal O}\/$ being standard.  There is a unique connection 
$\hbox{$ \nabla : {\cal F} $} \rightarrow {\cal F} \otimes \Omega\/$ compatible 
with the bilinear form on ${\cal F}\/$ such that the composition 
${\cal E} \hookrightarrow {\cal F} \righttypearrow{\nabla}{} 
{\cal F} \otimes \Omega \rightarrow {\cal E} \otimes \Omega\/$ is equal to
$\nabla_{\cal E}\/$ and $\nabla (1) = f\/$ where 
$1 \in H^0 (Y, {\cal O}) \subset H^0 (Y, {\cal F})\/$ and 
$f \in H^0 \bigl( Y, \Omega^{\otimes k} \bigr) = 
H^0 (Y, {\cal E}_1 \otimes \Omega) \subset 
H^0 (Y, {\cal F} \otimes \Omega)\/$.  Set ${\cal F}_i = {\cal E}_i\/$ for
$i \leq k-1\/$, ${\cal F}_i = {\cal E}_i \oplus {\cal O}\/$ for $i \leq k\/$.
The isomorphism ${\cal O} \righttypearrow{\sim}{} \det {\cal E}\/$ induces an
isomorphism ${\cal O} \righttypearrow{\sim}{} \det {\cal F}\/$.  So we have 
obtained a structure of \hbox{$S O (n)\/$-oper} on ${\cal F}\/$.  This 
\hbox{$S O (n)\/$-oper} will be denoted by ${\cal F} (L, f)\/$.  If ${\cal M}\/$ is a square
root of ${\cal O}\/$ $\bigl( \/$i.e., an invertible sheaf with a fixed isomorphism 
${\cal O} \righttypearrow {\sim}{} {\cal M}^{\otimes 2} \bigr)\/$ then 
${\cal F} (L, f) \otimes {\cal M}\/$ is also an \hbox{$S O (n)\/$-oper}.

\nullsubsection
\beginproclaim {Proposition}
Every \hbox{$S O (n)\/$-oper} is isomorphic to ${\cal F} (L, f) \otimes {\cal M}\/$ 
for some $L\/$, $f\/$ and ${\cal M}\/$.  ${\cal F} (L_1, f_1) \otimes 
{\cal M}_1\/$ and ${\cal F} (L_2, f_2) \otimes {\cal M}_2\/$ are isomorphic 
if and only if $L_1 = L_2\/$, $f_1 = f_2\/$ and there exists an isomorphism 
$\varphi : {\cal M}_1 \righttypearrow{\sim}{} {\cal M}_2\/$ such that the 
diagram 
\beginproeq
$$
\matrix { & {\cal O} & \cr
          \mapsw {\sim}{} \hskip -5.2em & & 
               \hskip -0.5em \mapse{\sim} \hfill  \cr
          {\cal M}_1 \otimes {\cal M}_1 & \hskip -1.0em 
               \maplongright {{\scriptstyle\varphi \otimes \varphi} \atop
                   {\displaystyle\sim}} 
                   \hskip -1.0em &
               {\cal M}_2 \otimes {\cal M}_2    \cr
        }
\eqnum {\twentytwo}{}
$$
\endproeq
is commutative.  Moreover, isomorphisms ${\cal F} (L, f) \otimes {\cal M}_1 
\righttypearrow{\sim}{} {\cal F} (L, f) \otimes {\cal M}_2\/$ bijectively 
correspond to isomorphisms $\varphi : {\cal M}_1 \righttypearrow{\sim}{} {\cal M}_2\/$ 
such that \label{\twentytwo}{} is commutative.
\endproclaim

\nullsubsection
\beginproclaim {Corollary}
There is a natural bijection between the set of \hbox{$s o (n)\/$-opers} and the set 
of pairs $(L, f)\/$ where $f \in H^0 \bigl( Y, \Omega^{\otimes k} \bigr)\/$ 
and $L : \Omega^{\otimes (1-k)} \rightarrow \Omega^{\otimes k}\/$ is a 
differential operator of order $2 k - 1\/$ such that $L^t = -L\/$ and the
principal symbol of $L\/$ equals~1.
\endproclaim

\nullsubsection
\beginremarks {Remark}
As explained in [DS1], [DS2], it is natural to associate to $(L, f)\/$ the
expression $L + f d^{-1} f \in H^0 
\bigl( Y, \Omega^{\otimes k} \otimes_{\cal O} \underline{\cal D} 
\otimes_{\cal O} \Omega^{\otimes (k-1)} \bigr)\/$ where $\underline{\cal D}\/$ 
is the sheaf of pseudodifferential symbols and 
$d^{-1} \in H^0 
\bigl( Y, \underline{\cal D} \otimes_{\cal O} \Omega^{\otimes (-1)} \bigr)\/$ 
is inverse to the differential operator $d : {\cal O} \rightarrow \Omega\/$ 
considered as an element of $H^0 
\bigl( Y, \Omega \otimes_{\cal O} {\cal D} \bigr)\/$.  One may say informally 
that $L + f d^{-1} f \in H^0 
\bigl( Y, \Omega^{\otimes k} \otimes_{\cal O} \underline{\cal D} 
\otimes_{\cal O} {\cal F}_1 \bigr)\/$ annihilates 
${\cal F}_1 \subset {\cal F}\/$.
\endremarks

The proof of the Proposition is based on the following lemma.

\nullsubsection
\beginproclaim {Lemma}
Let $({\cal F}, \{ {\cal F}_i \}, \nabla)\/$ be an \hbox{$S O (n)\/$-oper}.  Set 
${\cal M} = \Ker \bigl( \nabla : {\cal F}_k / {\cal F}_{k-1} \rightarrow 
({\cal F}_{k+1} / {\cal F}_k) \otimes \Omega \bigr)\/$.  Then 

\itemitem{1)} the bilinear form on ${\cal F}_k / {\cal F}_{k-1}\/$ induces 
              an isomorphism ${\cal M} \otimes {\cal M} \righttypearrow{\sim}{} {\cal O}\/$;

\itemitem{2)} {\rightskip = \parindent 
              there is a unique way of lifting the inclusion 
              ${\cal M} \hookrightarrow {\cal F}_k / {\cal F}_{k-1}\/$ to a 
              morphism $s : {\cal M} \rightarrow {\cal F}_k\/$ such that 
              $\nabla \bigl( s ({\cal M}) \bigr) \subset 
              \bigl( s ({\cal M}) + {\cal F}_1 \bigr) \otimes \Omega\/$. \par}

\endproclaim

Using the lemma it is easy to prove the proposition by considering the 
orthogonal decomposition ${\cal F} = s ({\cal M}) \oplus {\cal E}\/$.

\beginproof {\ of the lemma}
Since the composition of the morphisms ${\cal F}_{k-1} / {\cal F}_{k-2}
\rightarrow\break 
({\cal F}_k / {\cal F}_{k-1}) \otimes \Omega \rightarrow 
({\cal F}_{k+1} / {\cal F}_k) \otimes \Omega^{\otimes 2}\/$ induced by
$\nabla\/$ is an isomorphism ${\cal F}_k / {\cal F}_{k-1} = 
{\cal M} \oplus {\cal N}\/$ where ${\cal N} \otimes \Omega = \Im 
\bigl( \nabla : {\cal F}_{k-1} / {\cal F}_{k-2} \rightarrow 
({\cal F}_k / {\cal F}_{k-1}) \otimes \Omega \bigr)\/$.  Since $\nabla\/$
is compatible with the bilinear form ${\cal M}\/$ and ${\cal N}\/$ are 
orthogonal to each other.  Since the bilinear form on 
${\cal F}_k / {\cal F}_{k-1}\/$ is non-degenerate the statement 1) is clear.  
To prove the statement 2) we may assume that the curve $Y\/$ is affine and 
there is an isomorphism ${\cal O} \righttypearrow{\sim}{} {\cal M}\/$ 
compatible with the bilinear form on ${\cal M}\/$ (otherwise replace 
${\cal F}\/$ by ${\cal F} \otimes {\cal M}\/$).  Choose 
$\sigma \in H^0 (Y, {\cal F}_k)\/$ so that the image of $\sigma\/$ in
$H^0 (Y, {\cal F}_k / {\cal F}_{k-1})\/$ equals 
$1 \in H^0 (Y, {\cal O}) = H^0 (Y, {\cal M}) \subset 
H^0 (Y, {\cal F}_k / {\cal F}_{k-1})\/$.  We must show that there is a
unique $\tau \in H^0 (Y, {\cal F}_{k-1})\/$ such that 
$\nabla (\sigma + \tau) \in 
\bigl( {\cal O} (\sigma + \tau) \oplus {\cal F}_1 \bigr) \otimes \Omega\/$.
Since $(\sigma + \tau, \sigma + \tau) = 1\/$ we have 
$\bigl( \nabla (\sigma + \tau), \sigma + \tau \bigr) = 0\/$ and therefore if
$\nabla (\sigma + \tau) \in 
\bigl( {\cal O} (\sigma + \tau) \oplus {\cal F}_1 \bigr) \otimes \Omega\/$
then $\nabla (\sigma + \tau) \in {\cal F}_1 \otimes \Omega\/$.  Now the
uniqueness of $\tau\/$ is clear and the existence of $\tau\/$ means that
$$
\nabla (\sigma) \in \bigl( {\cal F}_{k-1} \otimes \Omega \bigr) +
\nabla \bigl( {\cal F}_{k-1} \bigr) 
\eqnum {\twentythree}{}
$$
Since the image of $\sigma\/$ in $H^0 (Y, {\cal F}_k)\/$ belongs to
${\cal M}\/$ we have $\nabla (\sigma) \in {\cal F}_k \otimes \Omega\/$ and since 
$\bigl( \nabla (\sigma), \sigma \bigr) = {1\over2} d (\sigma, \sigma) = 0\/$
the image of $\nabla (\sigma)\/$ in 
$({\cal F}_k / {\cal F}_{k-1}) \otimes \Omega\/$ belongs to 
${\cal N} \otimes {\cal M}\/$.  This is equivalent to \label{\twentythree}{}.
\endproof

\beginsection 3 A description of the set of \hbox{$\gg\/$-opers}

Denote by $\Op_\gg (Y)\/$ the set of isomorphism classes of \hbox{$\gg\/$-opers} 
on $Y\/$.

\subsection {3.1}
Let $\gg\/$= sl(2).  In this case it is well known that $\Op_\gg (Y)\/$
is a principal homogeneous space over 
$H^0 \bigl( Y, \Omega^{\otimes 2} \bigr)\/$.  Let us recall why it is so.  
According to \chapter.2.6 $\Op_\gg (Y)\/$ can be interpreted as the set of 
Sturm-Liouville operators $L\/$ on $Y\/$.  So we obtain a free transitive 
action of $H^0 \bigl( Y, \Omega^{\otimes 2} \bigr)\/$ on
$\Op_\gg (Y)\/$ defined by 
$$
L \mapsto L + \omega, \quad \omega \in H^0 \bigl( Y, \Omega^{\otimes 2} \bigr)
\eqnum {\twentyfour}{}
$$
Other interpretations of $\Op_\gg (Y)\/$ (see \chapter.2.6) give rise to 
actions of $H^0 \bigl( Y, \Omega^{\otimes 2} \bigr)\/$ on $\Op_\gg (Y)\/$ which 
can be written in terms of the Sturm-Liouville operators as 
$L \mapsto L + \alpha \omega\/$ where $\alpha \in {\Bbb Q}\/$ depends on the 
interpretation.

Let us recall why $\Op_\gg (Y) \neq \emptyset\/$.  Certainly there is an open 
covering $Y = \Cup_i U_i\/$ such that \hbox{$\gg\/$-opers} on ${\cal M}_i\/$ 
exist.  So if $H^1 \bigl( Y, \Omega^{\otimes 2} \bigr) = 0\/$ then 
$\Op_\gg (Y) \neq \phi\/$.  If 
$H^1 \bigl( Y, \Omega^{\otimes 2} \bigr) \neq 0\/$ then either $Y\/$ is an 
elliptic curve or $Y = {\Bbb P}^1\/$.  In the first case $\Omega\/$ is trivial 
and so the interpretation of \hbox{$\gg\/$-opers} as Sturm-Liouville operators 
shows that $\Op_\gg (Y) \neq \phi\/$.  On ${\Bbb P}^1\/$ there is an obvious 
\hbox{$S L (2)\/$-oper}: take the trivial bundle ${\cal E} = {\cal O}^2\/$ with 
the trivial connection and define ${\cal E}_1 \subset {\cal E}\/$ to be 
the tautological line subbundle.  So $\Op_\gg ({\Bbb P}^1) \neq \phi\/$.

\subsection {3.2}
The action \label{\twentyfour}{} can be described without using Sturm-Liouville 
operators.  Let $({\cal F}, \nabla)\/$ be a \hbox{$\gg\/$-oper} in the sense of
Definitions \chapter.1.2 and \chapter.1.3 where $\gg = s l (2)\/$.  Then we can
obtain a new \hbox{$\gg\/$-oper} $({\cal F}, \widetilde{\nabla})\/$ by setting 
$\widetilde{\nabla} = \nabla + \eta\/$, 
$\eta \in \Hom \bigl( \Theta, \gg^1_{\cal F} \bigr) = 
H^0 \bigl( Y, \gg^1_{\cal F} \otimes \Omega \bigr)\/$ where $\gg^1_{\cal F}\/$ 
has the 
same meaning as in \chapter.1.1.  We have the canonical isomorphism 
$\Theta \rightarrow \gg^{-1}_{\cal F} / \gg^0_{\cal F}\/$ induced by
$\nabla\/$ (see \chapter.1.2) and the pairing between 
$\gg^{-1}_{\cal F} / \gg^0_{\cal F}\/$ and $\gg^1_{\cal F}\/$ induced by 
the bilinear form.  $T r (X Y)\/$ on $\gg\/$.  So we can identify 
$\gg^1_{\cal F}\/$ with $\Omega\/$.  It is easy to see (cf. \chapter.1.4) that
the action \label{\twentyfour}{} corresponds to replacing $\nabla\/$ by 
$\nabla - \omega\/$, $\omega \in H^0 \bigl( Y, \Omega^{\otimes 2} \bigr) = 
H^0 \bigl( Y, \gg^1_{\cal F} \otimes \Omega \bigr)\/$.

\subsection {3.3}
Let $B \subset S L (2)\/$ denote the subgroup of upper-triangular matrices and 
$B_{\ad}\/$ the image of $B\/$ in $P G L (2)\/$.  It follows from \chapter.3.1 
and \chapter.3.2 that for any \hbox{$s l (2)\/$-opers} $({\cal F}, \nabla)\/$ and
$(\widetilde{\cal F}, \widetilde{\nabla})\/$ the \hbox{$B_{\ad}\/$-bundles} 
${\cal F}\/$ and $\widetilde{\cal F}\/$ are canonically isomorphic.  This
also follows from the construction of the \hbox{$S L (2)\/$-oper} corresponding 
to a Sturm-Liouville operator 
$\Omega^{\otimes (-1/2)} \rightarrow \Omega^{\otimes (-3/2)}\/$ 
(see \chapter.2.1 and \chapter.2.3 or \chapter.2.8).  Indeed, the
\hbox{$B\/$-bundle} in the construction corresponds to the 
\hbox{${\cal O}\/$-module} 
${\cal E} = {\cal D}_1 \otimes_{\cal O} \Omega^{\otimes (1/2)}\/$ and the 
canonical exact sequence $0 \rightarrow \Omega^{\otimes (1/2)} 
\rightarrow {\cal E} \rightarrow \Omega^{\otimes (-1/2)} \rightarrow 0\/$ 
(recall that ${\cal D}_1\/$ denotes the sheaf of differential operators of 
degree $\leq 1\/$).  This \hbox{$B\/$-bundle} depends on the choice of 
$\Omega^{\otimes (1/2)}\/$, but the corresponding \hbox{$B_{\ad}\/$-bundle} is 
a quite canonical object.

\nullsubsection
\beginremarks {Remark}
In [De1] Deligne gave a nice description of this \hbox{$B_{\ad}\/$-bundle} considered 
as a principal bundle:  its fiber over a point $y \in Y\/$ is the set of 
\hbox{${\Bbb C}\/$-isomorphisms} between $\Spec {\cal O}_y / m^3_y\/$ and 
$\Spec {\cal O}_\infty / m^3_\infty\/$ where ${\cal O}_y\/$ and 
${\cal O}_\infty\/$ are the local rings of $y \in Y\/$ and 
$\infty \in {\Bbb P}^1\/$ while $m_y\/$ and $m_\infty\/$ are their maximal
ideals ($B_{\ad}\/$ acts on this set because it acts on ${\Bbb P}^1\/$
preserving $\infty\/$).
\blackbox
\endremarks

\subsection {3.4}
Now let $\gg\/$ be an arbitrary semisimple Lie algebra with a fixed Borel 
subalgebra ${\frak b} \subset \gg\/$ and a fixed Cartan subalgebra 
${\frak h} \subset {\frak b}\/$.  Let $G_{\ad}\/$ be the connected algebraic 
group with trivial center corresponding to $\gg\/$.  Denote by $H_{\ad}\/$ and 
$B_{\ad}\/$ the Cartan and Borel subgroups of $G_{\ad}\/$ corresponding to 
${\frak h}\/$ and ${\frak b}\/$.  Fix an \hbox{$s l (2)\/$-triple} 
$\{ h,x,y \} \subset \gg\/$ such that $h \in {\frak h}\/$
and $x \in {\frak b}\/$ is a principle nilpotent element.  This means that 
$h = 2 \Sum_{\alpha \in \Gamma} \check{\omega}_\alpha\/$, 
$x = \Sum_{\alpha \in \Gamma} x_\alpha\/$, 
$y = \Sum_{\alpha \in \Gamma} y_\alpha\/$, where $\Gamma\/$ is the set of 
simple roots, $\check{\omega}_\alpha\/$ are the fundamental coweights, 
$x_\alpha \in \gg^\alpha\/$, $y_\alpha \in \gg^{-\alpha}\/$,
$x_\alpha \neq 0\/$, $y_\alpha \neq 0\/$.  The $y_\alpha\/$ are 
uniquely determined by $x_\alpha\/$ because $[x,y] = h\/$ while the 
elements $x_\alpha \in \gg^\alpha\/$, $x_\alpha \neq 0\/$, can be arbitrary.
However any two \hbox{$s l (2)\/$-triples} of the kind described above are 
conjugate by a unique element of $H_{\ad}\/$.

Set $V = \Ker \ad x \subset \gg\/$.  The grading 
$\gg = \relbuilder {\oplus}{k} \gg_k\/$ (see \chapter.1.1) induces a grading 
$V = \relbuilder {\oplus}{k} V_k\/$.  It is well known that $V_k = 0\/$ if
$k \leq 0\/$.  The numbers $k\/$ such that $V_k \neq 0\/$ are nothing but the 
exponents of $\gg\/$ and $\dim V\/$ is equal to the rank of $\gg\/$ 
(see [Ko1]).  Notice that $V_1 \neq 0\/$ because $x \in V_1\/$.  Of course
$V\/$ considered as a subspace of $\gg\/$ depends on $\{ h,x,y \}\/$.  But 
for any two choices of $\{ h,x,y \}\/$ there is a canonical isomorphism between 
the corresponding spaces $V\/$.  So we can identify all these spaces and 
consider the space $V (\gg)\/$ which depends only on $\gg\/$ (of course
$V (\gg)\/$ is not a subspace of $\gg\/$).  There is a canonical grading 
$V (\gg) = \relbuilder {\oplus}{k} V_k (\gg)\/$ and a canonical element 
${\cal X} \in V_1 (\gg)\/$ corresponding to $x \in V_1\/$.  Now set
$$
V^\gg_k (Y) = V_{k-1} (\gg) \otimes H^0 \bigl( Y, \Omega^{\otimes k} \bigr), 
\quad V^\gg (Y) = \relbuilder {\oplus}{k} V^\gg_k (Y)
\eqnum {\twentyfive}{}
$$
There is a canonical embedding $H^0 \bigl( Y, \Omega^{\otimes 2} \bigr) 
\hookrightarrow V^\gg_2 (Y) \hookrightarrow V^\gg (Y)\/$ defined by 
$$
\omega \mapsto -\omega \cdot {\cal X}
\eqnum {\twentysix}{}
$$
So the principal homogeneous space $\Op_{s l (2)} (Y)\/$ over 
$H^0 \bigl( Y, \Omega^{\otimes 2} \bigr)\/$ induces a principal homogeneous 
space $\underline{\Op}_\gg (Y)\/$ over $V^\gg (Y)\/$.  In other words 
$\underline{\Op}_\gg (Y)\/$ is the quotient of 
$\Op_{s l (2)} (Y) \times V^\gg (Y)\/$ modulo the equivalence relation
$$
(L + \omega, \eta) \sim (L, \eta - \omega {\cal X}) 
\eqnum {\twentyseven}{}
$$
where $\eta \in V^\gg (Y)\/$, 
$\omega \in H^0 \bigl( Y, \Omega^{\otimes 2} \bigr)\/$ 
and $L \in \Op_{s l (2)} (Y)\/$ is considered as a Sturm-Liouville operator.
We will construct a canonical bijection 
$\underline{\Op}_\gg (Y) \righttypearrow{\sim}{} {\Op}_\gg (Y)\/$.

\nullsubsection
\beginremarks {Remark}
Since $\Op_{\gg_1 \times \gg_2} (Y) = \Op_{\gg_1} (Y) \times \Op_{\gg_2} (Y)\/$
it is enough to describe $\Op_\gg (Y)\/$ if $\gg\/$ is simple.  In this case
$V_1 (\gg)\/$ is generated by ${\cal X}\/$ and therefore 
$\underline{\Op}_\gg (Y) = \Op_{s l (2)} (Y) \times 
\bigl( \relbuilder{\oplus}{k>2} V^\gg_k (Y) \bigr)\/$.
\blackbox
\endremarks

Let us construct a canonical mapping 
$\underline{\Op}_\gg (Y) \rightarrow \Op_\gg (Y)\/$.  By definition, an 
\hbox{$s l (2)\/$-oper} $L\/$ is a \hbox{$P G L (2)\/$-oper} 
$({\cal F}_0, \nabla_0)\/$ where ${\cal F}_0\/$ is a torsor over the 
upper-triangular subgroup $B_0 \subset P G L (2)\/$ and $\nabla_0\/$ is a 
connection on the corresponding \hbox{$P G L (2)\/$-bundle} satisfying a 
certain condition (see \chapter.1.2).  The fixed \hbox{$s l (2)\/$-triple} 
$\{ h,x,y \}\/$ defines a homomorphism $\rho : P G L (2) \rightarrow G_{\ad}\/$.  
Since $\rho (B_0) \subset B_{\ad}\/$ we can associate to ${\cal F}_0\/$ a 
\hbox{$B_{\ad}\/$-bundle} ${\cal F}\/$ and $\nabla_0\/$ induces a connection 
$\nabla\/$ on the \hbox{$G_{\ad}\/$-bundle} corresponding to ${\cal F}\/$.  
Clearly $({\cal F}, \nabla)\/$ is a \hbox{$G_{\ad}\/$-oper}, i.e., a 
\hbox{$\gg\/$-oper}.  For 
$\eta \in \Hom (\Theta, {\frak b}_{\cal F}) = 
H^0 (Y, {\frak b}_{\cal F} \otimes \Omega)\/$ 
the pair $({\cal F}, \nabla + \eta)\/$ is still a \hbox{$\gg\/$-oper}.  Notice 
that ${\frak b}_{\cal F} = {\frak b}_{{\cal F}_0}\/$ where 
${\frak b}_{{\cal F}_0}\/$ corresponds to 
${\cal F}_0\/$ and the representation of $B_0\/$ in ${\frak b}\/$ defined as the 
composition of $\rho : B_0 \rightarrow B_{\ad}\/$ and the adjoint representation 
of $B_{\ad}\/$ in~${\frak b}\/$.  Since $V\/$ is a \hbox{$B_0\/$-invariant} 
subspace of ${\frak b}\/$ we have 
${\frak b}_{\cal F} = {\frak b}_{{\cal F}_0} \supset V_{{\cal F}_0}\/$.  An 
element $g \in B_0\/$ acts on $V_k\/$ as multiplication by $\chi (g)^k\/$ 
where $\chi : B_0 \rightarrow {\Bbb G}_m\/$ maps $\minimatrix{a}{0}{b}{c}\/$ 
to $a/c\/$.  Since the line bundle corresponding to ${\cal F}_0\/$ and 
$\chi\/$ is canonically isomorphic to $\Omega\/$ (see \chapter.1.4) we have 
$V_{{\cal F}_0} = \relbuilder{\oplus}{k} 
\bigl( V_k \otimes \Omega^{\otimes k} \bigr)\/$.  So 
$H^0 \bigl( Y, {\frak b}_{\cal F} \otimes \Omega \bigr) \subset
H^0 \bigl( Y, V_{{\cal F}_0} \otimes \Omega \bigr) =
\relbuilder{\oplus}{k} \bigl( V_k \otimes H^0 
(Y, \Omega^{\otimes (k+1)}) \bigr) = V^\gg (Y)\/$.  Therefore to an 
\hbox{$s l (2)\/$-oper} $L\/$ and an $\eta \in V^\gg (Y)\/$ we can associate 
the \hbox{$\gg\/$-oper} $({\cal F}, \nabla + \eta)\/$ where 
$({\cal F}, \nabla)\/$ is
constructed from $L\/$ as explained above.  It follows from \chapter.3.2 and
\label{\twentyseven}{} that the constructed mapping 
$\Op_{s l (2)} (Y) \times V^\gg (Y) \rightarrow \Op_\gg (Y)\/$ factors through 
$\underline{\Op}_\gg (Y)\/$.  So we have obtained a canonical mapping
$$
\underline{\Op}_\gg (Y) \rightarrow \Op_\gg (Y)
\eqnum {\twentyeight}{}
$$

\nullsubsection
\beginproclaim {Theorem}
The mapping \label{\twentyeight}{} is bijective.
\endproclaim

\beginproof{}
Proposition \chapter.1.3 shows that the assertion is local.  So it is enough
to prove the following.  Suppose that there is function 
$z \in H^0 (Y, {\cal O}_Y)\/$ such that $dz\/$ has no zeros and let $q\/$
be a regular function on 
$Y\/$ with values in $\gg^{-1} = \relbuilder{\oplus}{k \geq -1} \gg_k\/$
such that for every simple root $\alpha\/$ the component $q_{- \alpha}\/$ of
$q\/$ has no zeros.  Then there is a unique $B_{\ad}\/$-valued regular function
$b\/$ on $Y\/$ such the values of the function $\widetilde{q}\/$ defined by
the relation 
$b^{-1} (\partial_z + q) b = \partial_z + \widetilde{q}\/$ belong to
$y + V\/$ where $y\/$ is the member of the \hbox{$s l (2)\/$-triple} 
$\{ h,x,y \}\/$ and $V = \Ker \ad x\/$.  This is easily deduced from the well
known equalities $\gg = (\Im \ad y) \oplus V\/$ and 
$(\Ker \ad y) \cap {\frak b} = 0\/$.
\endproof

Though the bijection \label{\twentyeight}{} is canonical we are not sure that it
gives a reasonable description of $\Op_\gg (Y)\/$.  So we will avoid using it
when possible.

\subsection {3.5}
Let $B_0\/$, $B_{\ad}\/$ and $\rho : P G L (2) \rightarrow G_{\ad}\/$ denote
the same objects as in \chapter.3.4.  Denote by ${\cal F}^{\can}_0\/$ the 
canonical \hbox{$B_0\/$-bundle} on $Y\/$ constructed in \chapter.3.3 and by 
${\cal F}^{\can}\/$ the \hbox{$B\/$-bundle} corresponding to ${\cal F}^{\can}_0\/$ and
the homomorphism $B_0 \rightarrow B\/$ induced by $\rho\/$.  It follows from
\chapter.3.3 and \chapter.3.4 that for any \hbox{$\gg\/$-oper} $({\cal F}, \nabla)\/$
the \hbox{$B_{\ad}\/$-bundle} ${\cal F}\/$ is canonically isomorphic to 
${\cal F}^{\can}\/$.

\subsection {3.6}
Theorem
\chapter.3.4 implies that for any finite subset $D \subset Y\/$ the natural 
mapping $\Op_\gg (Y) \rightarrow \Op_\gg (Y \backslash D)\/$ is injective.  We 
will usually consider $\Op_\gg (Y)\/$ as a subset of 
$\Op_\gg (Y \backslash D)\/$.

\nullsubsection
\beginremarks {Remark}
The injectivity of the mapping $\Op_\gg (Y) \rightarrow \Op_\gg 
(Y \backslash D)\/$ also follows from the definition of \hbox{$\gg\/$-oper} and the 
following fact:  if ${\cal E}_1\/$ and ${\cal E}_2\/$ are bundles with 
connections on $Y\/$ and $\varphi\/$ is a horizontal isomorphism between the 
restrictions of ${\cal E}_1\/$ and ${\cal E}_2\/$ to $Y \backslash D\/$ then 
$\varphi\/$ is the restriction of an isomorphism 
${\cal E}_1 \righttypearrow {\sim}{} {\cal E}_2\/$.  However the argument 
based on Theorem \chapter.3.4 works in more general situations: for 
     \hbox{$\gg\/$-opers} with singularities (see Proposition \chapter.4.2) and for 
\hbox{$(\gg, h)\/$-opers} (see \chapter.5.2).
\endremarks

\subsection {3.7}
The reader can skip this section and pass directly to \chapter\S 4.

As explained in \chapter.3.4 the space $V^\gg (Y)\/$ which appears in the 
definition of $\underline{\Op}_\gg (Y)\/$ is equal to 
$H^0 \bigl( Y, V_{{\cal F}_0} \bigr)\/$ where ${\cal F}_0\/$ is the canonical 
torsor over the upper-triangular subgroup $B_0 \subset P G L (2)\/$
(see \chapter.3.3) and $V = \Ker \ad x \subset \gg\/$ is considered as
a $B_0\/$ module via the embedding $B_0 \hookrightarrow B_{\ad}\/$ corresponding
to $\{ h,x,y \}\/$.  Now let $V'\/$ be any subspace of $\gg\/$ such that 
$$
[h, V'] \subset V', \quad
[x, V'] \subset V', \quad
\gg = (\Im \ad y) \oplus V', \quad
x \in V'
\eqnum {\twentynine}{}
$$
Then one can construct an analog of $\underline{\Op}_\gg (Y)\/$ in which 
$V^\gg (Y) = H^0 \bigl( Y, V_{{\cal F}_0} \bigr)\/$ is replaced by 
$H^0 \bigl( Y, {V'}_{{\cal F}_0} \bigr)\/$  and an analog of the bijection 
\label{\twentyeight}{}.  So it is natural to ask how to describe all 
subspaces $V' \subset \gg\/$ satisfying \label{\twentynine}{}.  If $V'\/$ is
such a subspace then we have the isomorphism 
$f : V' \righttypearrow {\sim}{} V\/$ defined as the composition 
$V' \righttypearrow {\sim}{} V / \Im \ad y \righttypearrow {\sim}{} V\/$.
Therefore $\ad x : V' \rightarrow V'\/$ induces the operator 
$\xi : V \rightarrow V\/$, $\xi = f \circ (\ad x) \circ f^{-1}\/$.  So
we obtain a mapping from the set of all $V'\/$ satisfying 
\label{\twentynine}{} to the set of all linear operators 
$\xi : V \rightarrow V\/$ such that 
$$
\xi (V_k) \subset V_{k + 1}, \quad 
\xi (x) = 0
\eqnum {\thirty}{}
$$
It is easy to show that this mapping is a bijection.  A glance into the list of 
exponents of simple Lie algebras shows that for simple Lie algebras $\gg\/$ of 
types $B_n\/$, $C_n\/$, $D_{2k}\/$, $E_7\/$, $E_8\/$, $F_4\/$, and $G_2\/$ the 
only possibility is $\xi = 0\/$, i.e., $V' = V\/$.  For $D_{2k + 1}\/$ and
$E_6\/$ there are some nontrivial choices of $V'\/$ and for $\gg = s l (n)\/$ 
there are plenty subspaces $V'\/$ satisfying \label{\twentynine}{}.  For 
instance one can take $V' = {\Bbb C} x \oplus {\Bbb C} e_{1n} \oplus \ldots 
\oplus {\Bbb C} e_{n-2, n}\/$ where $e_{ij}\/$ is the matrix with 1 at the 
intersection of the $i\/$-th row and the $j\/$-th column and zeros
elsewhere.

\beginsection 4 Singularities of \hbox{$\gg\/$-opers}

Just as in \chapter\S 3 we denote by $\gg\/$ a semisimple Lie algebra and by
$\Op_\gg (Y)\/$ the set of \hbox{$\gg\/$-opers} on $Y\/$.  For every finite subscheme 
$D \subset Y\/$ we will introduce a certain subset 
$\Op_{\gg, D} (Y) \subset \Op_\gg (Y \backslash D)\/$ such that 
$D \subset \Delta\/$ implies $\Op_{\gg, D} (Y) \subset \Op_{\gg, \Delta} (Y)\/$ 
and for every finite subset $S \subset Y\/$ \ $\Op_\gg (Y \backslash S)\/$ is 
the union of $\Op_{\gg, D} (Y)\/$ for all $D\/$ such that $\Supp D \subset S\/$.
Here $\Supp D\/$ denotes $D\/$ considered as a set.

\subsection {4.1}
In \chapter.3.4 we constructed a vector space $V^\gg (Y)\/$ 
(see \label{\twentyfive}{}) and a canonical embedding 
$H^0 (Y, \Omega^{\otimes 2}) \hookrightarrow V^\gg (Y)\/$ 
(see \label{\twentysix}{}).  We defined $\underline{\Op}_\gg (Y)\/$ to be the 
principal homogeneous space over $V^\gg (Y)\/$ induced by $\Op_{s l (2)} (Y)\/$
via this embedding $\bigl($recall that $\Op_{s l (2)} (Y)\/$ is a principal 
homogenous space over $H^0 (Y, \Omega^{\otimes 2}) \bigr)\/$.  Now for any
finite subscheme $D \subset Y\/$ set 
$$
V^\gg_{D,k} (Y) = V_{k-1} (\gg) \otimes H^0 
\bigl( Y, \Omega^{\otimes k} (k D) \bigr), \quad
V^\gg_D (Y) = \relbuilder {\oplus}{k} V^\gg_{D,k} (Y) 
\eqnum {\thirtyone}{}
$$
where $V_{k - 1} (\gg)\/$ has the same meaning as in \label{\twentyfive}{}.
Since $V^\gg (Y) \subset V^\gg_D (Y)\/$ the principal homogenous space
$\underline{\Op}_\gg (Y)\/$ over $V^\gg (Y)\/$ induces a principle homogenous 
space over $V^\gg_D (Y)\/$ which will be denoted by 
$\underline{\Op}_{\gg, D} (Y)\/$.  If $D \subset \Delta\/$ then 
$\underline{\Op}_{\gg, D} (Y) \subset \underline{\Op}_{\gg, \Delta} (Y)\/$ 
and for any finite subset $S \subset Y\/$ 
$\underline{\Op}_\gg (Y \backslash S)\/$ is the union of 
$\underline{\Op}_{\gg, D} (Y)\/$ 
for all $D\/$ such that $\Supp D \subset S\/$.

\subsection {4.2}
We have a canonical bijection 
$\underline{\Op}_\gg (Y \backslash D) \righttypearrow {\sim}{} \Op_\gg 
(Y \backslash D)\/$ (see \label{\twentyeight}{}).  So we could define 
$\Op_{\gg, D} (Y)\/$ as the image of 
$\underline{\Op}_{\gg, D} (Y) \subset \underline{\Op}_\gg (Y \backslash D)\/$
in $\Op_\gg (Y \backslash D)\/$.  But we prefer to give an equivalent 
``definition'' of $\Op_{\gg, D} (Y)\/$ which does not use the 
``suspicious'' bijection \label{\twentyeight}{}.

Let $G\/$ be the connected algebraic group with trivial center 
corresponding to $\gg, B\/$ the Borel subgroup of $G\/$ corresponding
to ${\frak b} \subset \gg\/$.  In the following definition we use the notation of
\chapter.1.1 and \chapter.1.2.

\nullsubsection
\beginproclaim {Definition}
A \hbox{$\gg\/$-oper} on $Y\/$ with $D\/$-singularities is a \hbox{$B\/$-bundle} ${\cal F}\/$
on $Y\/$ with a morphism 
$\nabla : \Theta (-D) \rightarrow {\cal E}^{-1}_{\cal F}\/$ such that the 
diagram
\beginproeq
$$
\matrix { \Theta (-D) & \hskip -1.0em \mapright{\nabla} \hskip -1.0em &
                        {\cal E}^{-1}_{\cal F}    \cr
          \hskip 2.0em \maphookdown & & \hskip -0.5em \mapsw {}{} \hfill  \cr
          \hfill \Theta \hskip -0.5em   \cr
        }
\eqnum {\thirtytwo}{}
$$
\endproeq
is commutative and for every simple root $\alpha\/$ the composition 
$\Theta (-D) \righttypearrow{\nabla}{} {\cal E}^{-1}_{\cal F} 
\righttypearrow{}{} {\cal E}^{-1}_{\cal F} / {\cal E}_{\cal F} =
\gg^{-1}_{\cal F} / \gg_{\cal F} 
\righttypearrow{}{} \gg^{-\alpha}_{\cal F}\/$ is an isomorphism.
\endproclaim

Denote by $\Op_{\gg, D} (Y)\/$ the set of isomorphism classes of
\hbox{$\gg\/$-opers} on $Y\/$ with\break 
\hbox{$D\/$-singularities.}

\nullsubsection
\beginproclaim {Proposition}
The natural mapping $\Op_{\gg, D} (Y) \rightarrow \Op_\gg (Y \backslash D)\/$
is injective and its image is equal to the image of 
$\underline{\Op}_{\gg, D} (Y) \subset \underline{\Op}_{\gg, D} (Y \backslash D)\/$
in ${\Op}_\gg (Y \backslash D)\/$.
\endproclaim

\beginproof{}
Fix an \hbox{$s l (2)\/$-triple} $\{ h,x,y \} \subset \gg\/$ such that 
$h \in {\frak h}\/$ and $x \in {\frak b}\/$ is a principle nilpotent element.  
Set $V = \Ker \ad x\/$.  We may assume that there is function 
$z \in H^0 (Y, {\cal O}_Y)\/$ such that $dz\/$ has no zeros and that 
$D \subset Y\/$ is defined by the equation $f = 0\/$, 
$f \in H^0 (Y, {\cal O}_Y)\/$.  The main step of the proof is to show that if
$q\/$ is a regular function on $Y\/$ with values in 
$\gg^{-1} = \relbuilder {\oplus}{k \geq -1} \gg_k\/$ such that for every
simple root $\alpha\/$ the component $q_{- \alpha}\/$ of $q\/$ has no
zeros then there is a unique $B\/$-valued regular function $b\/$ on $Y\/$ 
such that the values of the function $\widetilde{q}\/$ defined by the relation 
$b^{-1} (f \partial_z + q) b = f \partial_z + \widetilde{q}\/$ belong to
$y + V\/$.  This is proved just as in the case $f = 1\/$.  Then one has to find
the ``normal form'' of the operator $\partial_z + f^{-1} \widetilde{q}\/$,
i.e., one must find a $B\/$-valued regular function $g\/$ on
$Y \backslash D\/$ such that 
$g^{-1} (\partial_z + f^{-1} \widetilde{q}) g = \partial_z + \overline{q}\/$
where the values of $\overline{q}\/$ belong to $y + V\/$.  It is easy to see
that $g\/$ is the image of 
$$
\left( \matrix { f & \partial_z (f)  \cr
                 0 & 1      \cr
               } \right)
$$
under the homomorphism $P G L (2) \rightarrow G\/$ corresponding to 
$\{ h,x,y \}\/$.  Therefore the components $\overline{q}_k\/$ and 
$\widetilde{q}_k\/$ with respect to the grading
$\gg = \relbuilder {\oplus}{k} \gg_k\/$ are related by
$$
\eqalign{
         \overline{q}_k & = f^{-k -1} \widetilde{q}_k, \enspace k > 1    \cr
                 \noalign {\vskip 1.0ex}
         \overline{q}_1 & = f^{-2} \Bigl( q_1 +
               \bigl({\hbox{$\scriptstyle {1}\over{2}$}} 
                       \partial^2_z (f) \cdot f - 
                     {\hbox{$\scriptstyle {1}\over{4}$}}
                       \bigl( \partial_z (f) \bigr)^2 \bigr) z \Bigr)    \cr
         }
$$
So for any $k \geq 1\/$ the set of all possible $\overline{q}_k\/$ is
precisely 
$\gg_k \otimes H^0 \Bigl( Y, {\cal O}_Y \bigl( (k+1) D \bigr) \Bigr)\/$.
\endproof

Usually we will consider $\Op_{\gg, D} (Y)\/$ as a subset of 
$\Op_\gg (Y \backslash D)\/$.

\subsection {4.3}
Let $D\/$ and $D'\/$ be finite subschemes of $Y\/$ such that 
$D \subset D'\/$.  We will construct a mapping 
$\Op_{\gg, D} (Y) \rightarrow \Op_{\gg, D'} (Y)\/$ such that the diagram
$$
\matrix { \Op_{\gg, D} (Y) & \hskip -1.0em \mapright{} \hskip -1.0em &
                      \Op_{\gg, D'} (Y)     \cr
          \mapdown {\wr}  & &  \mapdown {\wr}      \cr
          \underline{\Op}_{\gg, D} (Y) & \hskip -1.0em \mapsubset
                             \hskip -1.0em &
                      \underline{\Op}_{\gg, D'} (Y)     \cr   
        }
\eqnum {\thirtythree}{}
$$
is commutative.

First of all to a \hbox{$B\/$-bundle} ${\cal F}\/$ and a finite subscheme
$\Delta \subset X\/$ we will associate a new \hbox{$B\/$-bundle} ${\cal F}^\Delta\/$. 
Let $U \subset Y\/$ be an open subset containing $\Delta\/$ such that the
restriction ${\cal F}_U\/$ of ${\cal F}\/$ to $U\/$ is trivial and there is
an $f \in H^0 (U, {\cal O}_U)\/$ such that $\Delta \subset U\/$ is defined by 
the equation $f = 0\/$.  Let $\lambda : {\Bbb G}_m \rightarrow H\/$ denote
the homomorphism such that for any simple root $\alpha\/$ \ $\lambda (t)\/$ 
acts on $\gg^\alpha\/$ as multiplication by $t\/$.  Fix a trivialization 
$\tau : {\cal F}_U \righttypearrow {\sim}{} U \times B\/$.  Let us explain
that by a \hbox{$B\/$-bundle} we understand a principle \hbox{$B\/$-bundle}, e.g., 
${\cal F}_U\/$ is a space with a free right action of $B\/$ such that 
${\cal F}_U / B = U\/$.  So the group of automorphisms of the \hbox{$B\/$-bundle}
$U \times B\/$ is the group of regular functions $U \rightarrow B\/$ which
acts on $U \times B\/$ by left multiplication.

\nullsubsection
\beginproclaim {Definition}
The \hbox{$B\/$-bundle} ${\cal F}^\Delta\/$ is obtained by gluing together 
${\cal F}_{Y \backslash \Delta}\/$ and $U \times B\/$ by means of the 
composition ${\cal F}_{U \backslash \Delta} \maplongright {\tau}
(U \backslash \Delta) \times B \maplongright {\lambda (f)} \vphantom {\Bigl(}
(U \backslash \Delta) \times B\/$.
\endproclaim

To check that ${\cal F}^{\Delta}\/$ is well defined we must show that if
$\tilde{\tau} : {\cal F}_U \righttypearrow {\sim}{} U \times B\/$ is 
another trivialization then the function 
$b : U \backslash \Delta \rightarrow B\/$ defined by 
$b = \bigl( \lambda (f) \tilde{\tau} \bigr) 
\bigl( \lambda (f) \tau \bigr)^{-1}\/$ is regular on~$U\/$.  
Indeed, $b = \lambda (f) (\tilde{\tau} \tau^{-1}) \lambda (f)^{-1}\/$,
$\tilde{\tau} \tau^{-1}\/$ is a regular function $U \rightarrow B\/$, and the 
coweight corresponding to $\lambda\/$ is dominant.

\nullsubsection
\beginremarks {Remark}
For a \hbox{$B\/$-module} $V\/$ we have a canonical grading 
$V = \relbuilder {\oplus}{k} V_k\/$ such that $\lambda (t) x = t^k x\/$ for
$x \in V_k\/$.  The filtration $V^k = \relbuilder {\oplus}{r \geq k} V_r\/$
is \hbox{$B\/$-invariant}.  So for any \hbox{$B\/$-bundle} 
${\cal F}\/$ on $Y\/$ it induces the filtration $V^k_{\cal F}\/$ of the 
vector bundle $V_{\cal F}\/$.  It is easy to show that 
$$
V_{{\cal F}^\Delta} = \Sum_k V^k_{\cal F} (k \Delta)
\eqnum {\thirtyfour}{}
$$
According to Theorem \chapter.3.2 from [DeM] a \hbox{$B\/$-bundle} ${\cal F}\/$ is
uniquely determined by the tensor functor $V \mapsto V_{\cal F}\/$ from the
category of finite dimensional \hbox{$B\/$-modules} to the category of vector
bundles on $Y\/$.  So \label{\thirtyfour}{} can be considered as another 
definition of ${\cal F}^\Delta\/$.
\blackbox
\endremarks

Let $({\cal F}, \nabla)\/$ be a \hbox{$\gg\/$-oper} on $Y\/$ with $D\/$-singularities
and $D'\/$ a finite subscheme on $Y\/$ containing $D\/$.  We are going to 
construct a \hbox{$\gg\/$-oper} $({\cal F}', \nabla')\/$ on $Y\/$ with 
$D'\/$-singularities.  Set ${\cal F}' = {\cal F}^\Delta\/$ where 
$\Delta \subset Y\/$ is a subscheme such that the corresponding ideal 
${\cal J}_\Delta \subset {\cal O}_Y\/$ is equal to ${\cal J}_{D'} (D)\/$ 
(in other words if $\Delta\/$, $D\/$, and $D'\/$ are considered as effective
divisors on $Y\/$ then $\Delta = D' - D\/$).  Let ${\cal E}_{\cal F}\/$ and
${\cal E}^{-1}_{\cal F}\/$ have the same meaning as in \chapter.1.1.  
${\cal F}'\/$ and ${\cal F}\/$ are canonically isomorphic outside $\Delta\/$.
So we have a canonical isomorphism between the restrictions of 
${\cal E}^{-1}_{\cal F}\/$ and ${\cal E}^{-1}_{{\cal F}'}\/$ to 
$Y \backslash D\/$.  It is easy to show that it is induced by a morphism 
${\cal E}^{-1}_{\cal F} \rightarrow {\cal E}^{-1}_{{\cal F}'} (\Delta)\/$ 
such that the corresponding mapping 
${\cal E}^{-1}_{\cal F} / {\cal E}_{\cal F} \rightarrow 
\bigl( {\cal E}^{-1}_{{\cal F}'} / {\cal E}_{{\cal F}'} \bigr) (\Delta)\/$ 
is an isomorphism.  We define $\nabla' : \Theta (-D') \rightarrow 
{\cal E}^{-1}_{{\cal F}'}\/$ to be the composition 
$\Theta (-D') \righttypearrow {\nabla}{} {\cal E}^{-1}_{\cal F} (- \Delta)
\rightarrow {\cal E}^{-1}_{{\cal F}'}\/$.

So for $D \subset D'\/$ we have constructed a canonical mapping 
$\Op_{\gg, D} (Y) \rightarrow \Op_{\gg, D'} (Y)\/$.  It is easy to see that
the diagram \label{\thirtythree}{} is commutative.  So the mapping 
$\Op_{\gg, D} (Y) \rightarrow \Op_{\gg, D'} (Y) \/$ is injective (this also
follows from Proposition \chapter.4.2).  We will usually identify 
$\Op_{\gg, D} (Y)\/$ with its image in $\Op_{\gg, D'} (Y)\/$.

\beginsection 5 The algebra $A_{\gg, D} (X)\/$

\subsection {5.1}
Let $X\/$ be a smooth projective curve over ${\Bbb C}\/$.  Then for any finite 
subscheme $D \subset X\/$ the affine space $\underline{\Op}_{\gg, D} (X)\/$ 
defined in \chapter.4.1 is finite dimensional.  So the bijection 
$$
\Op_{\gg, D} (X) \righttypearrow {\sim}{} \underline{\Op}_{\gg, D} (X)
\eqnum {\thirtyfive}{}
$$
defined in \chapter.4.2 induces a structure of algebraic variety on 
$\Op_{\gg, D} (X)\/$.  Of course this structure can be defined without 
using \label{\thirtyfive}{} because there is a natural notion of a family of 
\hbox{$\gg\/$-opers} with $D\/$-singularities parameterized by a scheme (the 
bijection \label{\thirtyfive}{} is just a convenient tool for proving that
the corresponding functor 
$\{${\it Schemes\/}$\}$ $\longrightarrow$ $\{${\it Sets\/}$\}$ is
representable by an algebraic variety isomorphic to an affine space).  On the
contrary, the structure of affine space on $\Op_{\gg, D} (X)\/$ induced by 
\label{\thirtyfive}{} is rather ``suspicious.''

Denote by $A_{\gg, D} (X)\/$ and $\underline{A}_{\gg, D} (X)\/$ the coordinate
rings of $\Op_{\gg, D} (X)\/$ and $\underline{\Op}_{\gg, D} (X)\/$.  We have 
the isomorphism
$$
A_{\gg, D} (X) \righttypearrow {\sim}{} \underline{A}_{\gg, D} (X)
\eqnum {\thirtysix}{}
$$
induced by \label{\thirtyfive}{}.  Of course $A_{\gg, D} (X)\/$ and 
$\underline{A}_{\gg, D} (X)\/$ are non-canonically isomorphic to 
${\Bbb C} [x_1, \ldots, x_N]\/$ where $N\/$ can be easily computed.

The vector space $V^\gg_D (X)\/$ corresponding to the affine space 
$\Op_{\gg, D} (X)\/$ has a canonical grading 
$V^\gg_D (X) = \relbuilder {\oplus}{k} V^\gg_{D, k} (X)\/$ 
(see \label{\thirtyone}{}).  Choosing a point 
$P \in \underline{\Op}_{\gg, D} (X)\/$ one obtains an isomorphism 
$\underline{\Op}_{\gg, D} (X) \righttypearrow {\sim}{} V^\gg_D (X)\/$ and
therefore a grading of $\underline{A}_{\gg, D} (X)\/$.  The 
increasing filtration $\underline{A}^{(n)}_{\gg, D} (X)\/$ corresponding to 
this grading does not depend on $P\/$.  It induces a filtration on
$A_{\gg, D} (X)\/$.  We will give a definition of this filtration which does 
not involve \label{\thirtyfive}{}.  We will also find the corresponding graded 
algebra.

\subsection {5.2}
Let $h \in {\Bbb C}\/$ (one can consider $h\/$ as ``Planck's constant'').  
Let $Y\/$ be a smooth curve over ${\Bbb C}\/$ and $G\/$ an algebraic group
over ${\Bbb C}\/$.  If ${\cal F}\/$ is a \hbox{$G\/$-bundle} and ${\cal E_F}\/$ 
the algebroid of infinitesimal symmetries of 
${\cal F}\/$ then by an $h\/$-{\it connection\/} on ${\cal F}\/$ we mean a 
morphism $\nabla : \Theta \rightarrow {\cal E_F}\/$ such that the diagram
$$
\matrix { \hskip 0.5em \Theta & \hskip -1.0em \mapright{h} \hskip -1.0 em & 
                 \Theta     \cr
          \hskip 0.5em \maptypedown {\nabla}{} & & 
                 \hskip -0.5em \mapnw {}{} \hfill       \cr
          \hfill {\cal E_F} \hskip -0.5em     \cr
        }
$$
is commutative.  If ${\cal L}\/$ is a locally free sheaf of 
\hbox{${\cal O}_Y\/$-modules} then an $h\/$-connection on ${\cal L}\/$ is a
${\Bbb C}\/$-linear morphism 
$\nabla : {\cal L} \rightarrow {\cal L} \otimes \Omega\/$ such that 
$\nabla (f s) = f \nabla (s) + h s \otimes d f\/$ where $f\/$ and $s\/$
are local sections of ${\cal O}\/$ and ${\cal L}\/$.

If $h = 1\/$ then an $h\/$-connection is just a connection.  The case 
$h \neq 0\/$ can be reduced to $h = 1\/$ (replace $\nabla\/$ by 
$h^{-1} \nabla\/$).  If $h = 0\/$ then an $h\/$-connection on a vector
bundle ${\cal L}\/$ is an ${\cal O}_Y\/$-linear morphism 
${\cal L} \rightarrow {\cal L} \otimes \Omega\/$ and an $h\/$-connection on a
principle \hbox{$G\/$-bundle} ${\cal F}\/$ is an element of 
$H^0 (Y, \gg_{\cal F} \otimes \Omega)\/$ where $\gg_{\cal F}\/$ is the vector 
bundle corresponding to ${\cal F}\/$ and the adjoint representation of $G\/$.

Now let $G\/$ be a connected reductive group over ${\Bbb C}\/$.  Replacing
the word ``connection'' in Definition \chapter.1.2 by ``$h\/$-connection'' 
one obtains the notion of \hbox{$(G, h)\/$-{\it oper\/}}.  In the same way one 
introduces the notions of \hbox{$(\gg, h)\/$-{\it oper\/}} and 
\hbox{$(\gg, h)\/$-{\it oper\/}} {\it with \hbox{$D\/$-singularities}.\/}  
Denote by $\Op^h_\gg (Y)\/$ $\bigl($resp. $\Op^h_{\gg, D} (Y) \bigr)\/$ the set 
of \hbox{$(\gg, h)\/$-opers} on $Y\/$ (resp. \hbox{$(\gg, h)\/$-opers} on
$Y\/$ with \hbox{$D\/$-singularities).}

Most of the results of \chapter\S 1--\chapter\S 4 concerning 
\hbox{$G\/$-opers}, \hbox{$\gg\/$-opers} and \hbox{$\gg\/$-opers} 
with\break
\hbox{$D\/$-singularities} 
are also true for \hbox{$(G,h)\/$-opers}, \hbox{$(\gg,h)\/$-opers}, 
and \hbox{$(\gg,h)\/$-opers} with\break
\hbox{$D\/$-singularities}  
(the only exception is the 
interpretation of analytic \hbox{$s l (2)\/$-opers} as locally projective 
structures in \chapter.2.7).  Of course one should replace connections by 
$h\/$-connections and differential operators by 
$h\/$-{\it differential operators\/}.

By definition, the sheaf of $h\/$-differential operators ${\cal D}_h\/$ is
the sheaf of rings generated by the ring ${\cal O}\/$ and the left 
\hbox{${\cal O}\/$-module} $\Theta\/$  with the following defining relations:
if $v \in \Theta\/$ (i.e., $v\/$ is a local section of $\Theta\/$) and 
$f \in {\cal O}\/$ then $v f = f v + h v (f)\/$ where $v (f)\/$ is the 
derivative of $f\/$ with respect to $v\/$ (if $\dim Y\/$ were greater than 1 
we would also have to add the relations 
$v_1 v_2 - v_2 v_1 = h [v_1, v_2]\/$ for $v_1, v_2 \in \Theta\/$ but if
$\dim Y = 1\/$ these relations hold automatically).  If $h = 1\/$ then 
${\cal D}_h = {\cal D}\/$ while if $h = 0\/$ then ${\cal D}_h\/$ is the
symmetric algebra of $\Theta\/$.  For any $h\/$ there is a canonical morphism
${\cal D}_h \rightarrow {\cal D}\/$ which is identical on ${\cal O}\/$ and
acts on $\Theta\/$ as multiplication by~$h\/$.  It is an isomorphism for
$h \neq 0\/$.

In \chapter.2.6 we associated to each \hbox{$s l (2)\/$-oper} a differential operator 
$\Omega^{\otimes (-1/2)} \rightarrow \Omega^{\otimes (3/2)}\/$ of order 2.
In the same way one associates to an \hbox{$(s l (2), h)\/$-oper} an 
$h\/$-differential operator 
$\Omega^{\otimes (-1/2)} \rightarrow \Omega^{\otimes (3/2)}\/$.  Applying
the canonical morphism ${\cal D}_h \rightarrow {\cal D}\/$ one obtains a 
differential operator 
$\Omega^{\otimes (-1/2)} \rightarrow \Omega^{\otimes (3/2)}\/$  with principal
symbol $h^2\/$.  Thus we get a natural bijection between 
$\Op^h_{s l (2)} (Y)\/$ and the set of differential operators 
$L : \Omega^{\otimes (-1/2)} \rightarrow \Omega^{\otimes (3/2)}\/$ of order
2 such that $L^t = L\/$ and the principal symbol of $L\/$ equals~$h^2\/$.
In particular $\Op^h_{s l (2)} (Y)\/$ is a principal homogenous space over
$H^0 \bigl( Y, \Omega^{\otimes 2} \bigr)\/$.  The description of the action of 
$H^0 \bigl( Y, \Omega^{\otimes 2} \bigr)\/$ on $\Op_{s l (2)} (Y)\/$ given in
\chapter.3.2 works also for $\Op^h_{s l (2)} (Y)\/$.

Now let $\gg\/$ be an arbitrary semisimple Lie algebra over ${\Bbb C}\/$ and
$V^\gg_D (Y)\/$ the vector space defined by \label{\thirtyone}{}.  Denote by 
$\underline{\Op}^h_{\gg, D} (Y)\/$ the principal homogeneous space over 
$V^\gg_D (Y)\/$ induced by $\Op^h_{s l (2)} (Y)\/$ via the embedding 
$H^0 \bigl( Y, \Omega^{\otimes 2} \bigr) \hookrightarrow V^\gg_D (Y)\/$ 
defined by \label{\twentysix}{}.  Just as in \chapter.3.4 and \chapter.4.2
one constructs a canonical bijection
$$
\underline{\Op}^h_{\gg, D} (Y) \righttypearrow {\sim}{}
\Op^h_{\gg, D} (Y) 
\eqnum {\thirtyseven}{}
$$

\subsection {5.3}
Let $X\/$ be a smooth projective curve over ${\Bbb C}\/$ and $D\/$ a 
finite subscheme of $X\/$.  Set 
$$
\matrix{
       \hfill \Opscript_{\gg, D} (X) & \hskip -0.5em = \hskip -0.5em &
                    \Bigl\{ ( h,u ) \bigm| h \in {\Bbb C}, \hfill & 
                    \hskip -0.5em u \in 
                          \Op^h_{\gg, D} (X) \Bigr\} \hfill      \cr
       \noalign {\vskip 1.0ex}
       \hfill \underline\Opscript_{\gg, D} (X) & 
                    \hskip -0.5em = \hskip -0.5em & 
                    \Bigl\{ ( h,u ) \bigm| h \in {\Bbb C}, \hfill & 
                    \hskip -0.5em u \in 
                          \underline{\Op}^h_{\gg, D} (X) \Bigr\} \hfill   \cr
       }
$$
There is a natural structure of algebraic variety on 
$\Opscript_{\gg, D} (X)\/$ and $\underline\Opscript_{\gg, D} (X)\/$ and 
\label{\thirtyseven}{} induces an isomorphism 
$\underline\Opscript_{\gg, D} (X) \righttypearrow {\sim}{} 
\Opscript_{\gg, D} (X)\/$.  So $\Opscript_{\gg, D} (X)\/$ is an affine
variety non-canonically isomorphic to an affine space and the morphism 
$\pi : \Opscript_{\gg, D} (X) \rightarrow {\Bbb A}^1\/$ which maps 
$\Op^h_{\gg, D} (X)\/$ to $h \in {\Bbb C}\/$ is flat.  ${\Bbb G}_m\/$ acts on
$\Opscript_{\gg, D} (X)\/$ (replacing $\nabla\/$ by $\lambda \nabla\/$,
$\lambda \in {\Bbb C}^*\/$) and $\pi\/$ is ${\Bbb G}_m\/$-equivariant.  Set
$\Opscript^*_{\gg, D} (X) = \pi^{-1} ({\Bbb G}_m)\/$.  We have 
$\Op_{\gg, D} (X) = \pi^{-1} (1) = 
{\Bbb G}_m \backslash \Opscript^*_{\gg, D} (X)\/$.  Therefore 
$A_{\gg, D} (X)\/$, the coordinate ring of $\Op_{\gg, D} (X)\/$, can be
identified with the ring of \hbox{${\Bbb G}_m\/$-invariant} regular functions on 
$\Opscript^*_{\gg, D} (X)\/$.  So we have a filtration 
$A_{\gg, D} (X) = \Cup_n A^{(n)}_{\gg, D} (X)\/$ where 
$A^{(n)}_{\gg, D} (X)\/$ is the space of \hbox{${\Bbb G}_m\/$-invariant} regular
functions on $\Opscript^*_{\gg, D} (X)\/$ having a pole of order $\leq n\/$
at the divisor $\pi^{-1} (0)\/$.  It is easy to show that this filtration
corresponds to the filtration 
$\underline{A}_{\gg, D} (X) = \Cup_n \underline{A}^{(n)}_{\gg, D} (X)\/$
defined in \chapter.5.1.

\subsection {5.4}
The graded ring $g r A_{\gg, D} (X) = 
\relbuilder {\oplus}{n} A^{(n)}_{\gg, D} (X) / A^{(n-1)}_{\gg, D} (X)\/$ is
nothing but the coordinate ring of $\pi^{-1} (0) = \Op^0_{\gg, D} (X)\/$.
The grading on the coordinate ring comes from the action of ${\Bbb G}_m\/$
on $\Op^0_{\gg, D} (X)\/$.  We will describe $\Op^0_{\gg, D} (X)\/$ 
together with this action.  In what follows we use the notation of 
\chapter.1.1 and \chapter.1.2.

Let $G\/$ be the connected algebraic group with trivial center corresponding 
to $\gg\/$, $B\/$ its Borel subgroup.  By definition, points of 
$\Op^0_{\gg, D} (X)\/$ correspond to pairs $({\cal F}, \omega)\/$ where 
${\cal F}\/$ is a \hbox{$B\/$-bundle} on $X\/$ and $\omega \in 
H^0 \bigl( X, \gg^{-1}_{\cal F} \otimes \Omega (D) \bigr)\/$ is such that for 
every simple root $\alpha\/$ the image of $\omega\/$ in 
$H^0 \bigl( X, \gg^{-\alpha}_{\cal F} \otimes \Omega (D) \bigr)\/$ has no 
zeros as a section of $\gg^{-\alpha}_{\cal F} \otimes \Omega (D)\/$.
Denote by $\Inv (\gg)\/$ the algebra of \hbox{$G\/$-invariant} polynomials on 
$\gg\/$.  If $f \in \Inv (\gg)\/$ is homogeneous of degree $n\/$ then
$f (\omega) \in H^0 \bigl( X, \Omega^{\otimes n} (n D) \bigr)\/$.  So we 
obtain a mapping 
$$
\Op^0_{\gg, D} (X) \rightarrow \Hitch_{\gg, D} (X) =
      \Mor (\Inv (\gg), \relbuilder {\oplus}{n}
                         H^0 \bigl( X, \Omega^{\otimes n} (n D) \bigr)
\eqnum {\thirtyeight}{}
$$
where $\Mor\/$ denotes the set of graded algebra morphisms.  There is a 
canonical isomorphism
$$
\Inv (\gg) \righttypearrow {\sim}{} \Inv (\gg^{\wedge *}) 
\eqnum {\thirtynine}{}
$$
where $\Inv (\gg^{\wedge *})\/$ is the algebra of invariant polynomials
of the Langlands dual $\gg^{\wedge}\/$ $\bigl($indeed, both 
$\Inv (\gg)\/$ and $\Inv (\gg^{\wedge *})\/$ are canonically isomorphic to the algebra of
polynomials on ${\frak h}\/$ invariant with respect to the Weyl group$\bigr)$.
So the r.h.s. of \label{\thirtyeight}{} can be identified with
$\Hitch_{\gg^{\wedge}, D} (X)\/$, the base of Hitchin's fibration
(see 0.2), and \label{\thirtyeight}{} can be interpreted 
as a mapping
$$
\Op^0_{\gg, D} (X) \rightarrow \Hitch_{\gg^\wedge, D} (X)
\eqnum {\forty}{}
$$
$\Hitch_{\gg^\wedge, D} (X)\/$ has a natural structure of affine algebraic
variety with an action of ${\Bbb G}_m\/$ on it (see section 0).

\nullsubsection
\beginproclaim {Proposition}
The mapping \label{\forty}{} is a ${\Bbb G}_m\/$-equivariant isomorphism of 
algebraic varieties.
\endproclaim

\nullsubsection
\beginproclaim {Corollary}
$g r A_{\gg, D} (X)\/$ is canonically isomorphic to the coordinate ring of
$\Hitch_{\gg, D} (X)\/$.
\endproclaim

The proof of the proposition is easy:  it is clear that \label{\forty}{} is a
${\Bbb G}_m\/$-equivariant morphism and the assertion that it is an isomorphism
is an immediate consequence of the following theorem proved in [Ko2].

\nullsubsection
\beginproclaim {Kostant's Theorem}
Set $\gg^{-1}_* = \Bigl\{ u + \Sum_{\alpha \in \Gamma} v_\alpha \bigm| 
u \in {\frak b}\/$, $v_\alpha \in \gg^{- \alpha}\/$, 
$v_\alpha \neq 0 \Bigr\}\/$ where $\Gamma\/$
is the set of simple roots.  Denote by $f\/$ the canonical mapping 
$\gg \rightarrow \Spec \Inv (\gg) = W \backslash {\frak h}\/$ where $W\/$ 
is the Weyl group.  Then

\itemitem{(i)} {\rightskip = \parindent 
               there is an affine subspace $L \subset \gg^{-1}_*\/$ such that
               the mapping $B \times L \rightarrow \gg^{-1}_*\/$, 
               $(b, v) \mapsto b v b^{-1}\/$, is an isomorphism of algebraic 
               varieties
               \par}

\itemitem{(ii)} {\rightskip = \parindent 
                the mapping $B \backslash \gg^{-1}_* = 
                L \rightarrow W \backslash {\frak h}\/$ induced by $f\/$ is an
                isomorphism of algebraic varieties.
                \par}

\endproclaim

For the reader's convenience we sketch a proof of Kostant's theorem.  Let 
$\{ h,x,y \}\/$ be the \hbox{$s l (2)\/$-triple} considered in \chapter.3.4.  Set 
$V = \Ker \ad x\/$ and $L = y + V\/$.  The proof of the statement (i) is 
quite similar to that of Proposition \chapter.1.3 and Theorem \chapter.3.4
(one uses the equalities $\gg = (\Im \ad y) \oplus V\/$ and 
$(\Ker \ad y) \cap {\frak b} = 0)\/$.  The mapping 
${\frak h} \rightarrow B \backslash \gg^{-1}_*\/$ which associates to 
$a \in {\frak h}\/$ the $b\/$-conjugacy class of $a + y\/$ factors through
$W \backslash {\frak h}\/$ because for every simple root $\alpha\/$
$$
s_\alpha (a) + y = \exp \bigl( - \alpha (a) \cdot \ad u_\alpha \bigr) 
                   (a + y), \quad
                   a \in {\frak h}
$$
where $s_\alpha (a) = a - \alpha (a) \check{\alpha}\/$ and 
$u_\alpha \in \gg^{\alpha}\/$ is such that $[u_\alpha, y] = \check{\alpha}\/$.
The composition $W \backslash {\frak h} \rightarrow 
B \backslash \gg^{-1}_* \rightarrow W \backslash {\frak h}\/$ is identical 
because the closure of the conjugacy class of $a + y\/$ contains $a\/$ for 
$a \in {\frak h}\/$.  It remains to show that 
$\dim (B \backslash \gg^{-1}_*) = \dim (W \backslash {\frak h})\/$,
i.e., $\dim V = \dim {\frak h}\/$.  Indeed, $\dim V = \dim \Ker \ad x = n\/$
where $n\/$ is the number of irreducible components of $\gg\/$ considered 
as a module over $s l (2) = {\Bbb C} h + {\Bbb C} x + {\Bbb C} y\/$.  But 
since the eigenvalues of $\ad h\/$ are even $n = \dim \Ker \ad h = \dim h\/$.

\nullsubsection
\beginremarks {Remark}
There is a canonical \hbox{${\Bbb G}_m\/$-invariant} isomorphism of algebraic varieties
$$
V^\gg_D (X) \righttypearrow {\sim}{} \Hitch_{\gg^\wedge, D} (X)
\eqnum {\fortyone}{}
$$
where $V^\gg_D (X)\/$ is defined by \label{\thirtyone}{} and 
$\lambda \in {\Bbb G}_m\/$ 
acts on $V^\gg_{D, k} (X)\/$ as multiplication by $\lambda^k\/$.  To construct
\label{\fortyone}{} choose a principal \hbox{$s l (2)\/$-triple} $\{ h,x,y \}\/$ and 
identify the spaces $V_{k - 1} (\gg)\/$ from \label{\thirtyone}{}) with the 
eigenspaces of $\ad h\/$ in $\Ker \ad x\/$.  Then to every 
$\eta \in V^\gg_D (X)\/$ and every homogeneous $f \in \Inv (\gg)\/$ of degree 
$n\/$ we associate 
$f (y + \eta) \in H^0 \bigl( X, \Omega^{\otimes n} (n D) \bigr)\/$.  This
defines the mapping \label{\fortyone}{} and Kostant's theorem shows that it is an 
isomorphism.
\endremarks

\bigskip\smallskip
\parskip=4pt
\baselineskip=13pt plus 2pt minus 2pt
\parindent=33pt
\frenchspacing
\noindent
{\larger\bf References}

\item{[De1]\hfill} P. Deligne, {\it \'{E}quations differentielles \`{a} points 
             singuliers r\'{e}guliers\/}, Lecture Notes in Mathematics, 
             Springer-Verlag, Berlin, vol.~{\bf 163} (1970).

\item{[DeM]\hfill} P. Deligne and J. Milne,  {\it Tannakian categories\/},  
             Lecture Notes in Mathematics, 
             Springer-Verlag, Berlin, vol.~{\bf 900} (1982), p.~101--228.

\item{[DS1]\hfill} V. G. Drinfeld and V. V. Sokolov, 
             {\it Equations of Korteweg-deVries type and simple Lie algebras\/},
             Soviet Mathematics Doklady, vol.~{\bf 23} (1981), No.~3, 
             p.~457--462.

\item{[DS2]\hfill} V. G. Drinfeld and V. V. Sokolov, 
             {\it Lie algebras and equations of Korteweg-deVries type\/},
             Journal of Soviet Mathematics, vol.~{\bf 30} (1985), p.~1975--2035.

\item{[Hi1]\hfill} N. J. Hitchin, {\it Stable bundles and integrable systems\/},
             Duke Math. Journal, vol.~{\bf 54} (1987), p.~91--114.

\item{[Ko1]\hfill} B. Kostant, {\it The principal three-dimensional subgroup and the
             Betti numbers of a complex simple Lie group\/},
             American Journal of Math., vol.~{\bf 81} (1959), p.~973--1032.

\item{[Ko2]\hfill} B. Kostant, {\it Lie group representations on polynomial rings\/},
             American Journal of Math., vol.~{\bf 85} No.~3, (1963), p.~327--404.

\item{[Te1]\hfill} C. Teleman, {\it Sur les structures homographiques d'une surface 
             de Riemann\/},  Comment. Math. Helv., vol.~{\bf 33} (1959), 
             p.~206--211.

\item{[Tyu]\hfill} A. N. Tyurin, {\it On the periods of quadratic differentials\/},
             Uspekhi Mat. Nauk, vol.~{\bf 33} No.~6, (1978), p.~149--195 
             (Russian).

\bye